\documentclass[fleqn,a4paper,10pt]{amsart}

\usepackage{graphicx}
\usepackage{tikz}
\usetikzlibrary{shapes,decorations.pathmorphing}
\usepackage{amsmath}
\usepackage{amssymb,amsthm}
\usepackage{comment}
\usepackage{hyperref}

\theoremstyle{definition}

\newtheorem{proper}{Property}

\title[Two-point function of triangulations:
a new derivation]{The distance-dependent two-point function of triangulations:
a new derivation from old results}
\author{Emmanuel Guitter}
\address{Institut de Physique Th\'eorique, CEA, IPhT, 91191 Gif-sur-Yvette, France, CNRS, UMR 3681}
\email{emmanuel.guitter@cea.fr}

\begin{document}
\maketitle

\begin{abstract}
We present a new derivation of the distance-dependent two-point function of random planar triangulations. 
As it is well-known, this function is intimately related to the generating functions of so-called slices, which are pieces of
triangulation having boundaries made of shortest paths of prescribed length. We show that the slice generating functions
are fully determined by a direct recursive relation on their boundary length. Remarkably, the kernel of this 
recursion is some quantity introduced and computed by Tutte a long time ago in the context of a global enumeration of planar 
triangulations. We may thus rely on these old results to solve our new recursion relation explicitly in a constructive way.
\end{abstract}

\section{Introduction}
\label{sec:introduction}

The combinatorics of planar maps, i.e.\ connected graphs embedded on the sphere, is making constant progress since the seminal 
work of Tutte in the $60$'th. In the more recent years, a growing interest was shown for metric properties of maps endowed with their 
graph distance, and especially for the corresponding distance statistics within ensembles of random maps. An emblematic result was 
the computation, for several families of maps, of the \emph{distance-dependent two-point function} which, so to say, measures the profile of
distances between two points (vertices or edges) picked at random on the map. Explicit expressions for this two-point function were obtained for
ensembles of planar maps with controlled face degrees \cite{GEOD,BG12} as well as for maps (or hyper-maps) with arbitrary face degrees
but with controlled edge (or hyper-edge) and face numbers \cite{AmBudd,BFG}. A first way to solve these questions was the use of bijections between maps and decorated trees, as first discovered by Schaeffer \cite{SchPhD} (upon reformulating a bijection by Cori and Vauquelin \cite{CoriVa}). In a second, intimately related, approach, 
the problem of computing the distance-dependent two-point function was reduced to that of enumerating \emph{slices}, which are particular pieces of 
maps bordered by shortest paths of prescribed length, meeting at some ``apex". In a first stage, the computation of either decorated trees or slice generating functions relied 
on finding the solutions of particular \emph{integrable systems of equations} satisfied by the generating functions at hand. 
No general technique however was developed to solve these equations and all the explicit expressions obtained in this way 
were the result of a \emph{simple guessing} of the solution. 
The recourse to decorated trees or slices took on its full dimension when it was later discovered that their generating functions could be obtained mechanically as coefficients in suitable continued fraction expansions for standard map generating functions.
This property was exploited in \cite{BG12} to obtain a constructive derivation of the distance dependent two-point function for maps with controlled
face degrees.

In this paper, we revisit the problem of computing the \emph{distance-dependent two-point function of random planar triangulations}, i.e\ planar
maps whose all faces have degree $3$. These maps were extensively studied in the past as they form one of the simplest natural families of maps.
Their two-point function was first obtained in \cite{PDFRaman} by guessing the solution of the associated integrable system. It was then
re-obtained in \cite{BG12} as a particular example of the general continued fraction formalism. Here, we present a \emph{new recursive
approach} which consists in directly relating the generating function of slices whose border has (maximal) length $k$ to that of
slices whose border has (maximal) length $k-1$ (see eq.~\eqref{eq:newrec} below). Remarkably enough, the ``kernel" of our recursion relation is some particular generating function of triangulations, already introduced by Tutte as early as in is first paper \cite{TutteCPT} on triangulations. We may thus directly use the old results 
of \cite{TutteCPT}  to solve our new recursion relation in a constructive way, without recourse to any guessing. 
 
The paper is organized as follows: in Section~\ref{sec:slicegf}, we recall the definition of slices and their connection with the distance-dependent
two-point function of random planar triangulations. We also recall the standard integrable system obeyed by the slice generating functions,
whose solution was guessed in \cite{PDFRaman}. Section~\ref{sec:newrec} is devoted to the derivation of our new recursion relation
between the generating function $T_k$ for slices with (maximal) border length $k$ and $T_{k-1}$. This new recursion is
based on the existence of some particular \emph{dividing line} which, so to say, delimits in the slice a region whose vertices are at distance
strictly larger than $k-1$ from the apex of the slice. As just mentioned, the kernel of our recursion is some particular 
generating function computed by Tutte in its seminal paper \cite{TutteCPT} on triangulations. In order to stick to Tutte's original results, 
we make in Section~\ref{sec:simpletriang} a detour to the family of \emph{simple triangulations}, i.e.\ triangulations with neither loops nor multiple edges. As shown, a simple substitution procedure makes the correspondence between this simplified family and the family of all triangulations that we are interested in.
We then use in Section~\ref{sec:Tuttesol} the explicit form given by Tutte for the kernel of our recursion relation to rewrite this recursion 
in a particularly simple and classical form (see eq.~\eqref{eq:Ykrec} below), whose solution is easily obtained by classical techniques.  
We finally return in Section~\ref{sec:finalexpr}
to the case of general (not necessarily simple) triangulations by performing the required substitution. This leads us to our final explicit expressions
for $T_k$ and for the distance-dependent two-point function.
We gather our concluding remarks in Section~\ref{sec:conclusion}.   

\section{Slice generating functions: reminders}
\label{sec:slicegf}
\subsection{Definitions}
\label{sec:slicedefs}
The distance-dependent two-point function of planar triangulations may be expressed in terms of the
generating functions $R_k$ and $S_k$ for \emph{$R$-slices} and \emph{$S$-slices} of maximal size $k$
(see eq.~\eqref{eq:twopoint} below). Slices are particular families of \emph{triangulations 
with a boundary}, namely planar rooted (i.e.\ with 
a marked oriented edge, the root-edge) maps whose all faces have degree $3$, except for the outer face (i.e.\ the face lying on the right of the root-edge) which may have arbitrary degree. The inner faces form what it called the \emph{bulk} while the edges incident to the outer face (visited, say clockwise around the bulk) form \emph{the boundary} whose \emph{length} is the degree of the outer face. $R$- and $S$-slices
are defined as follows:
\begin{figure}
\begin{center}
\includegraphics[width=12cm]{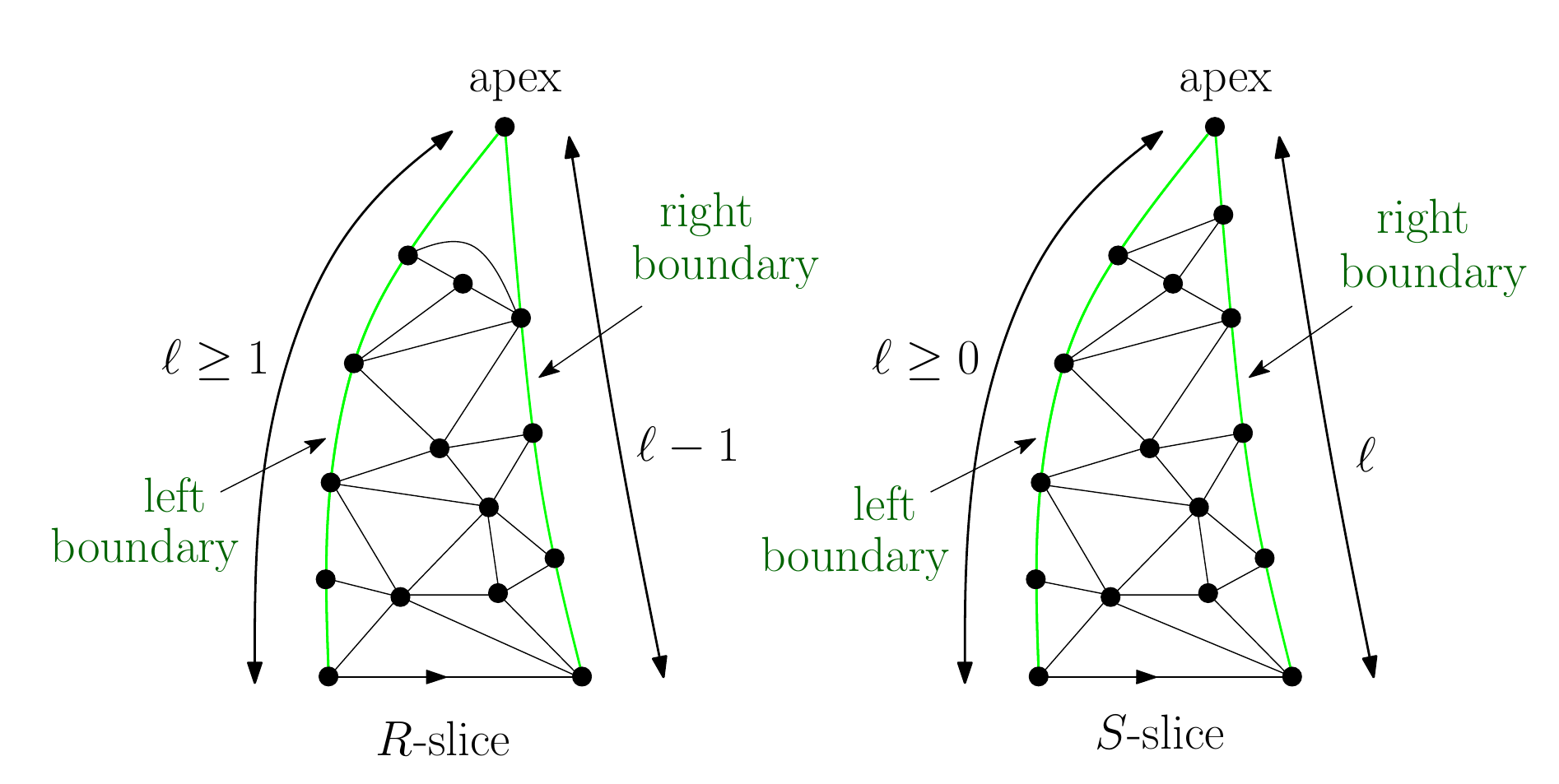}
\end{center}
\caption{An example of $R$-slice (left) and $S$-slice (right). In both cases, the left boundary is a shortest path within the map
between the root-vertex and the apex and has length $\ell$. As for the right boundary, it is in both cases the unique shortest 
path within the map between the endpoint of the root-edge and the apex, with length $\ell-1$ for the $R$-slice ($\ell\geq 1$), and $\ell$ for the $S$-slice ($\ell\geq 0$).}
\label{fig:slice}
\end{figure}

\begin{itemize}
\item $R$-slices have a boundary of length $2\ell$ ($\ell\geq 1$) and satisfy (see figure~\ref{fig:slice}):
 \begin{itemize}
 \item the (graph) distance from the origin of the root-edge (the root-vertex) to the \emph{apex}, which is  the vertex 
 reached from the root-vertex by making $\ell$ elementary steps  
 \emph{along the boundary} clockwise around the bulk, is $\ell$.
 In other words, the left boundary of the slice, which is the part of its boundary lying between the root-vertex and the apex clockwise around the bulk
 is a shortest path between its endpoints within the map;
 \item the distance from the endpoint of the root-edge to the apex is $\ell-1$.
 In other words, the right boundary of the slice, which is the part of the boundary lying between the endpoint of the root-edge and the apex 
 counterclockwise around the bulk is a shortest path between its endpoints within the map; 
 \item the right boundary is the \emph{unique} shortest path between its endpoints within the map;
 \item the left and right boundaries do not meet before reaching the apex.
 \end{itemize}
We call $R_k\equiv R_k(g)$ ($k\geq 1$) the generating function of $R$-slices with $1\leq \ell \leq k$, enumerated with a weight $g$
per inner face.
\vskip .2cm

\item $S$-slices have a boundary of length $2\ell+1$ ($\ell\geq 0$) and satisfy (see figure~\ref{fig:slice}):
 \begin{itemize}
 \item the distance from the root-vertex to the apex,  which is  the vertex 
 reached from the root-vertex by making $\ell$ elementary steps  
 \emph{along the boundary} clockwise around the bulk, is $\ell$.
 In other words, the left boundary of the slice (which is the part of the boundary lying between the root-vertex and the apex clockwise
 around the bulk) is a shortest path between its endpoints within the map;
 \item the distance from the endpoint of the root-edge to the apex is $\ell$.
 In other words, the right boundary of the slice (which is the part of the boundary lying between the endpoint of the root-edge and the apex 
 counterclockwise around the bulk) is a shortest path between its endpoints within the map; 
 \item the right boundary is the unique shortest path between its endpoints within the map;
 \item the left and right boundaries do not meet before reaching the apex.
 \end{itemize}
We call $S_k\equiv S_k(g)$ ($k\geq 0$) the generating function of $S$-slices with $0\leq \ell \leq k$, enumerated with a weight $g$
per inner face.
\end{itemize}
Note that the map reduced to a single root-edge and an outer face of degree  $2$ is an $R$-slice 
with $\ell=1$ and contributes a term $1$ to $R_k$ for any $k\geq 1$.
\vskip .2cm

\begin{figure}
\begin{center}
\includegraphics[width=9cm]{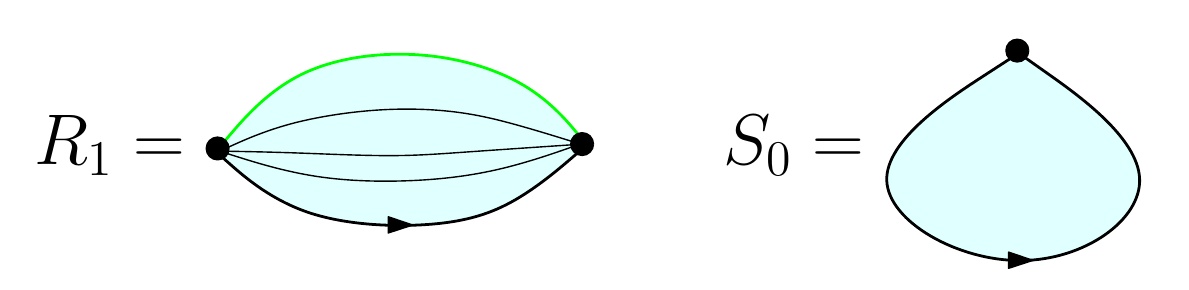}
\end{center}
\caption{A schematic picture of maps enumerated by $R_1$ and $S_0$.}
\label{fig:R1S0}
\end{figure}
Of particular interest are the generating functions $R_1$ and $S_0$, with the following interpretations: by definition, $R_1$ enumerates
$R$-slices with $\ell=1$, therefore with a boundary of length $2$. The right boundary has length $0$, hence the apex is
the endpoint of the root edge. The left boundary, of length $1$, connects the extremities of the root-edge, which are 
necessarily distinct. The function $R_1$ may therefore be interpreted as the \emph{generating function of
rooted triangulations with a boundary of length $2$ connecting two distinct vertices} (the extremities of the root-edge). 
Note that in such maps, the edges connecting the extremities of the root
edge within the map form in general what we shall call a a \emph{bundle} of $p$ edges for some $p\geq 1$ (see figure~\ref{fig:R1S0}). 
The function $R_1$ therefore enumerates bundles of edges. 
As for $S_0$, it enumerates
$S$-slices with $\ell=0$, in which case both extremities of the root-edge coincide with the apex. In particular,
the root-edge forms a loop. The function $S_0$ may thus be interpreted as the \emph{generating function of
rooted triangulations with a boundary of length $1$} (see figure~\ref{fig:R1S0}).

For later use, we also introduce the generating function
\begin{equation*}
T_k\equiv S_k-S_0
\end{equation*}
for $k\geq 0$. This generating function enumerates $S$-slices with a left-boundary length $\ell$ satisfying $1\leq \ell\leq k$. 
These are precisely the $S$-slices contributing to $S_k$ and \emph{whose root-edge does not form a loop} (recall indeed that the
left and right boundary of an $S$-slice are required to meet only at the apex, so that the root-edge forms a loop if and only if
$\ell=0$). Note that $T_0=0$ by definition.
\vskip .2cm

\begin{figure}
\begin{center}
\includegraphics[width=11cm]{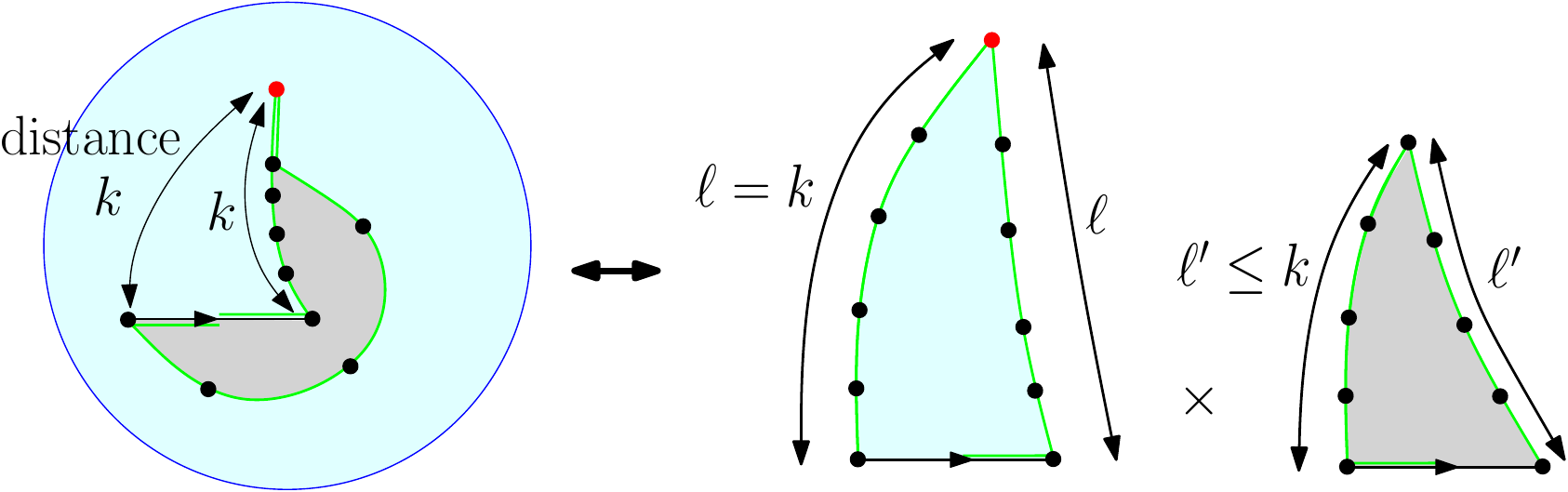}
\end{center}
\caption{A schematic picture of the one-to-one correspondence between pointed rooted triangulations whose root-edge 
has its both extremities at distance $k$ from the marked vertex (in red) and a pair of $S$-slices with left-boundary lengths $\ell$
and $\ell'$ satisfying $\max(\ell,\ell')=k$. The green lines on the left side are the leftmost shortest paths to the marked vertex
from the middle of the root-edge in both directions. Cutting along these paths creates the two slices on the right. The choice of
leftmost shortest paths ensures that the right boundaries of the slices are the unique shortest paths between their endpoints
within the slice. Note that we reversed the original orientation of the root-edge to obtain the root-edge of the grey slice.}
\label{fig:twopoint2}
\end{figure}
\begin{figure}
\begin{center}
\includegraphics[width=8cm]{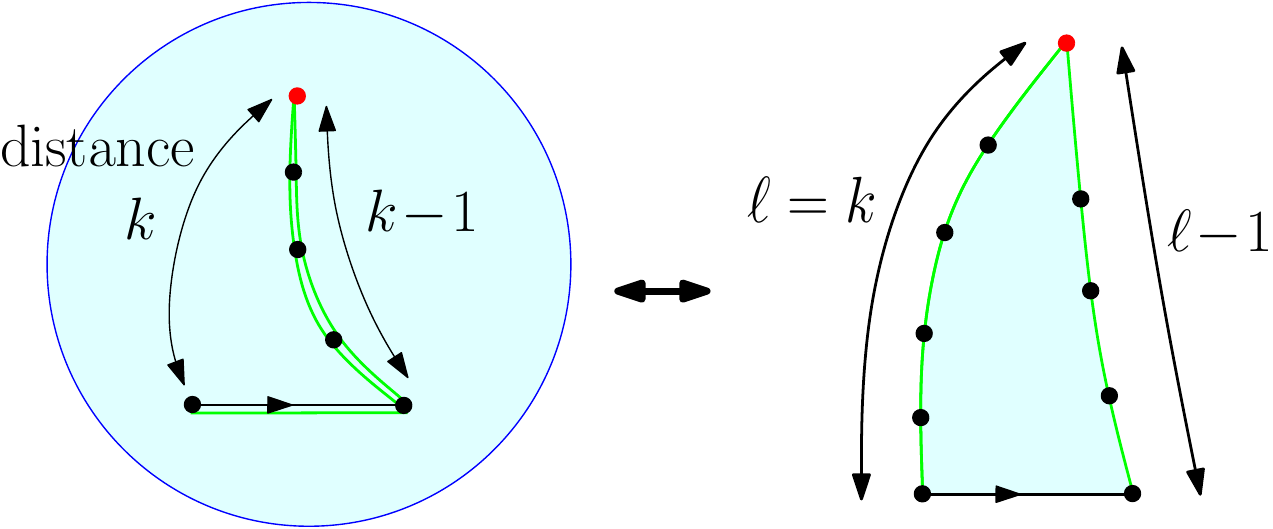}
\end{center}
\caption{A schematic picture of the one-to-one correspondence between a pointed rooted triangulation whose root-edge 
has its origin and endpoint at respective distance $k$ and $k-1$ from the marked vertex (in red) and an $R$-slices with left-boundary lengths $\ell=k$.
The green line on the left side is the leftmost shortest paths to the marked vertex
from the the root-vertex (with the root-edge as fist step). Cutting along this path creates the desired $R$-slice on the right.}
\label{fig:twopoint1}
\end{figure}
\subsection{The distance dependent two-point function} 
\label{sec:twopoint}
For $k\geq 0$, we define the distance-dependent two-point function of planar triangulations as the generating function $G_k$ 
for pointed (i.e.\ with a marked vertex) rooted (i.e.\ with a marked oriented edge) planar triangulations (i.e\ planar
maps whose all faces have degree $3$) for which the marked vertex is at graph distance $k$ from the
root-vertex (i.e.\ the origin of the root-edge). Let us show the following identity for $k\geq 1$:
\begin{equation}
G_k=(S_k^2-S_{k-1}^2)+(R_k-R_{k-1}-\delta_{k,1})+(R_{k+1}-R_k)=S_k^2-S_{k-1}^2+R_{k+1}-R_{k-1}-\delta_{k,1}
\label{eq:twopoint}
\end{equation}
with the convention $R_0=0$.
Maps enumerated by $G_k$ may indeed be classified into three classes according to the distance from their
marked vertex to the endpoint of the root-edge. This distance may be $k$, $k-1$ or $k+1$ and the three terms in the 
middle expression in \eqref{eq:twopoint} above correspond to the enumeration of the three classes. If the two extremities of the 
root-edge are at distance $k$ from the marked vertex,  we draw the \emph{leftmost} shortest paths to the marked vertex,
starting from the middle of the root-edge in both directions (see figure~\ref{fig:twopoint2}). Cutting along these paths results into 
two $S$-slices, whose root-edge is the original root-edge with its original orientation for one slice and
with the reversed orientation for the other. Note that the choice of leftmost shortest paths ensures that the right boundary of each piece 
is the unique shortest path between its endpoints within the piece. As for the left boundaries, they are also shortest paths between their endpoints,  with lengths less than or equal to $k$ (corresponding to situations enumerated by $S_k^2$)  and with at least one
of these lengths equal to $k$, hence are enumerated by $S_k^2-S_{k-1}^2$ (since the situations where both lengths are strictly
less than $\ell$ are enumerated by $S_{k-1}^2$).  If the endpoint of the
root-edge is at distance $k-1$ from the marked vertex,  we draw the leftmost shortest path form the root-vertex to the marked vertex,
taking the root-edge as first step (see figure~\ref{fig:twopoint1}). Cutting along this path results into 
an $R$-slice whose root-edge is the original root-edge, with left boundary length equal to $k$, and moreover different from the single root-edge if $k=1$.
Such slices are enumerated by $(R_k-1)-(R_{k-1}-1)=R_k-R_{k-1}$ for $k\geq 2$ and by $R_1-1$ for $k=1$,
i.e.\ by $R_k-R_{k-1}-\delta_{k,1}$ with our convention that $R_0=0$.
Finally, if the endpoint of the
root-edge is at distance $k+1$ from the marked vertex,  we reverse the orientation of the root-edge to get back to the previous
situation with $k\to k+1$. Such maps are thus enumerated by $R_{k+1}-R_k$, hence the relation \eqref{eq:twopoint}. 
As already mentioned, the knowledge of both $R_k$ and $S_k$ leads immediately via \eqref{eq:twopoint} to a expression for the 
distance-dependent two-point function $G_k$.
As a final remark, for $k=0$, $G_0$ enumerates rooted triangulations and adapting the above argument immediately yields
\begin{equation*}
G_0=S_0^2+R_1-1\ .
\end{equation*}
\vskip .2cm
\begin{figure}
\begin{center}
\includegraphics[width=11cm]{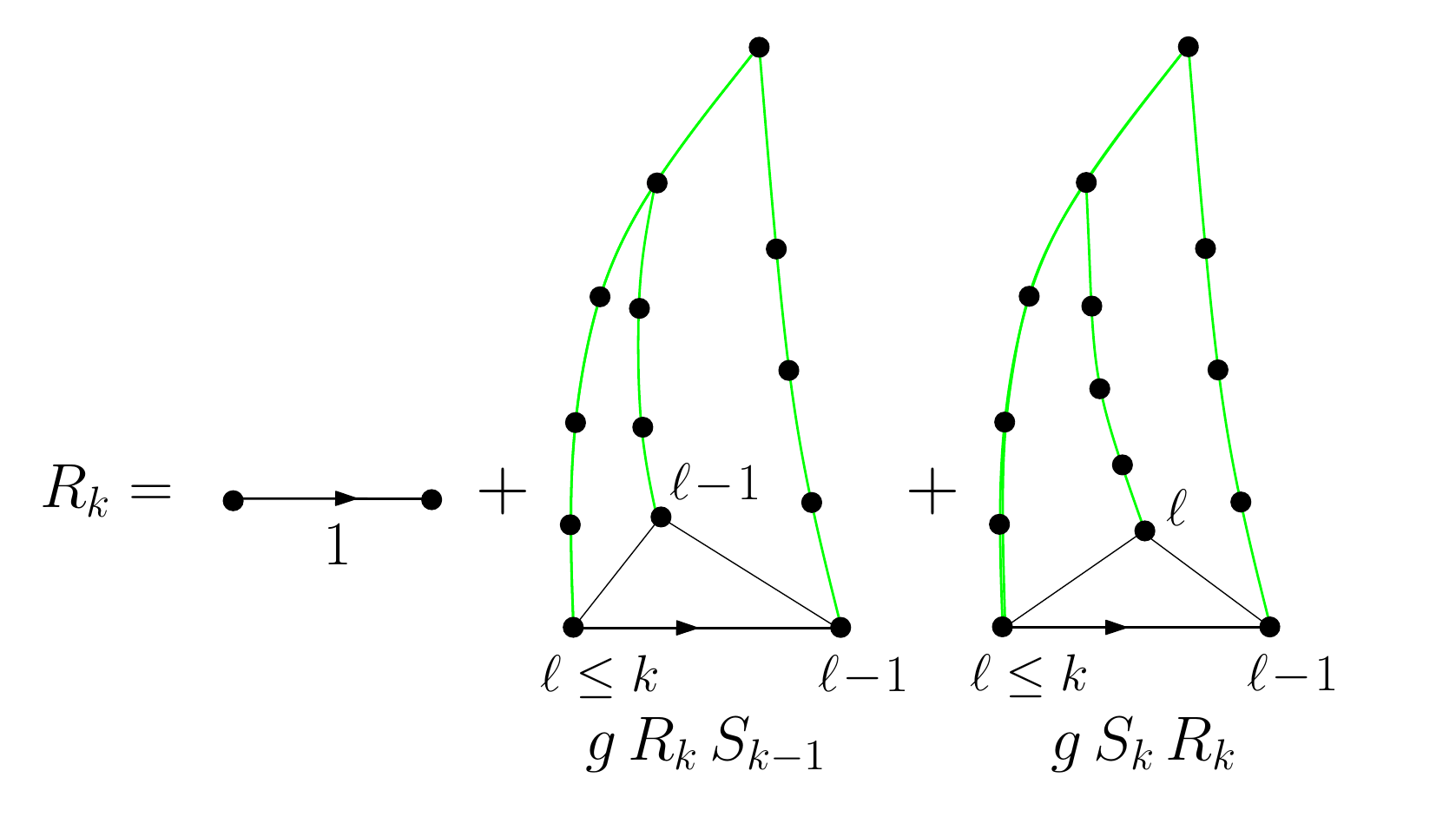}
\end{center}
\caption{A schematic picture explaining the first line of eq.~\eqref{eq:oldrec}. If not reduced to a single edge, the $R$-slice with $1\leq\ell\leq k$ is decomposed
by removing the triangle immediately on the left of the root-edge, whose intermediate vertex is at distance $\ell-1$ or $\ell$ from the apex , and
by cutting along the leftmost shortest path from this vertex to the apex.}
\label{fig:recurslice}
\end{figure}
\begin{figure}
\begin{center}
\includegraphics[width=10cm]{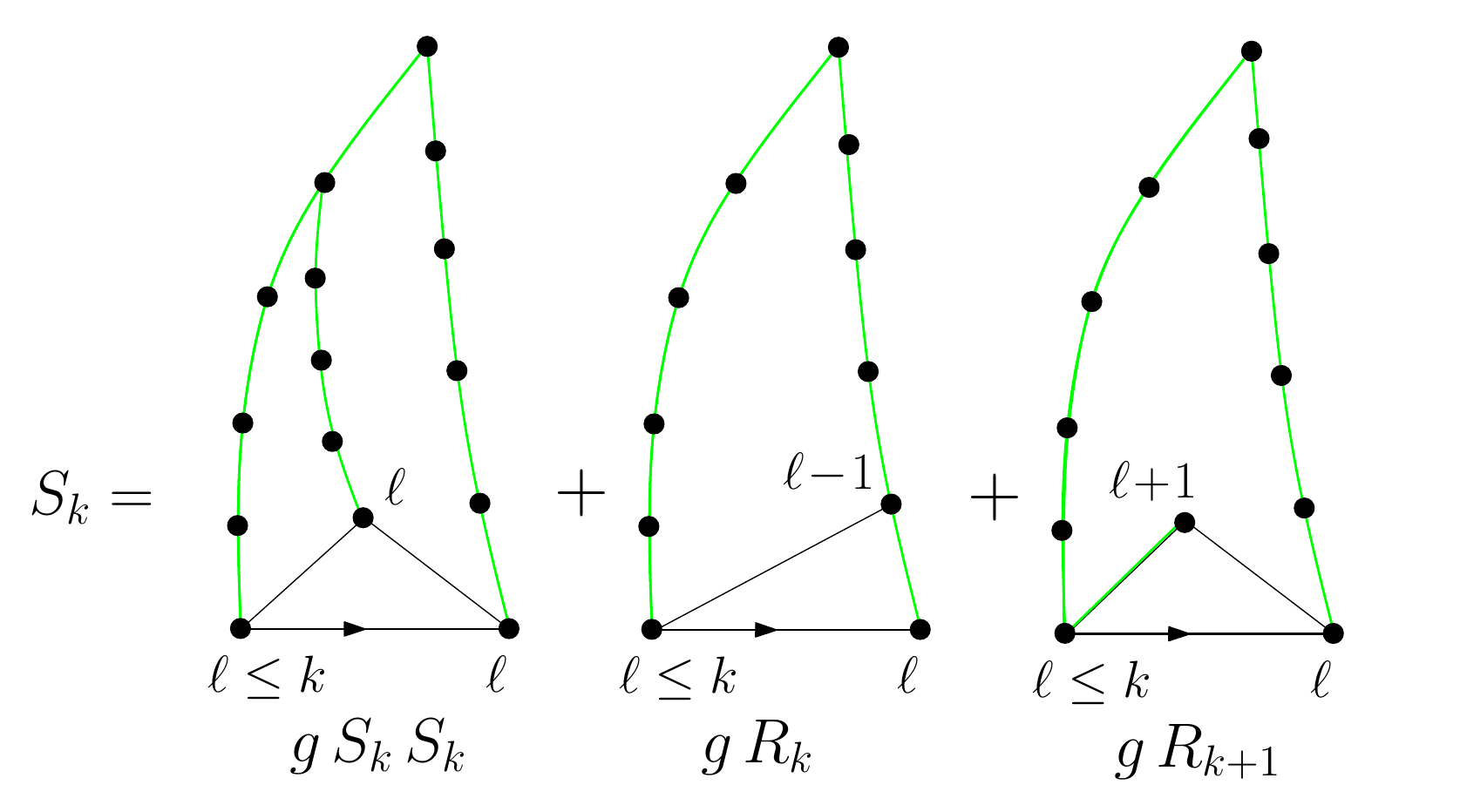}
\end{center}
\caption{A schematic picture explaining the second line of eq.~\eqref{eq:oldrec}. The $S$-slice with $0\leq\ell\leq k$ is decomposed
by removing the triangle immediately on the left of the root-edge, whose intermediate vertex is at distance $\ell$ or $\ell-1$ 
(in which case it lies on the right boundary) or $\ell+1$ from the apex , and, in the first case, 
by cutting along the leftmost shortest path from this vertex to the apex.}
\label{fig:recurslicebis}
\end{figure}

\subsection{Classical relations for slice generating functions} 
\label{sec:oldrec}
A first set of relations for the slice generating functions may be obtained by classifying the slices according to the nature of the inner face lying immediately on the left of the root-edge.
In the case of an $R$-slice not reduced to a single root-edge,  this triangular face is incident to the two 
extremities of the root-edge, at respective distances $\ell$ and $\ell-1$ from the apex, and to an intermediate vertex 
at distance $\ell-1$ or $\ell$ (note that this vertex may possibly be identical to one of the two others). Drawing
the leftmost shortest path from this intermediate vertex to the apex separates the map into two slices: an $R$-slice
and an $S$-slice (see figure \ref{fig:recurslice}). As for $S$-slices,  the triangular face on the left of the root-edge
has its intermediate vertex 
at distance $\ell$ or $\ell\pm 1$ from the apex. In the latter case, removing the triangle directly results into an $R$-slice
while, in the former case, drawing
the leftmost shortest path from the intermediate vertex to the apex separates the map into two $S$-slices
(see figure \ref{fig:recurslicebis}).
Taking into account the boundary length constraints, we immediately arrive at the classical system
\begin{equation}
\left\{
\begin{split}
R_k& =1+g\, R_k(S_{k-1}+S_k)\ , \qquad k\geq 1\\
S_k& = g\, (S_k^2 +R_k+R_{k+1})\ , \qquad k\geq 0\\
\end{split}
\label{eq:oldrec}
\right.
\end{equation}
with again our convention $R_0=0$. Note that this system is not sricto sensu recursive since we don't know the value 
of $S_0$. Still it is recursive order by order in $g$ if we impose $R_k=1+O(g^2)$ for all $k\geq 1$ and 
$S_k=O(g)$ for all $k\geq 0$, as required by the definition of $R_k$ and $S_k$ as slice generating functions.
The solution of this system was first found in \cite{PDFRaman} by simple guessing, leading to explicit expressions 
for $R_k$ and $S_k$, hence for the distance dependent two-point function $G_k$ via eq.~\eqref{eq:twopoint}. Later on, these explicit
expressions were re-derived in a constructive way {\it without recourse to the system \eqref{eq:oldrec}}
or to any other recursion relation, 
but by instead relating $R_k$ and $S_k$ to the distance-independent generating functions of
triangulations with a boundary of fixed length \cite{BG12}. 

In the next section, we shall introduce a new set of 
recursion relations for $R_k$ and $S_k$, which we shall then solve explicitly in a constructive way. 
 
\section{A new approach by recursion}
\label{sec:newrec}
\subsection{Construction of a dividing line}
\label{sec:dividing}
\begin{figure}
\begin{center}
\includegraphics[width=7cm]{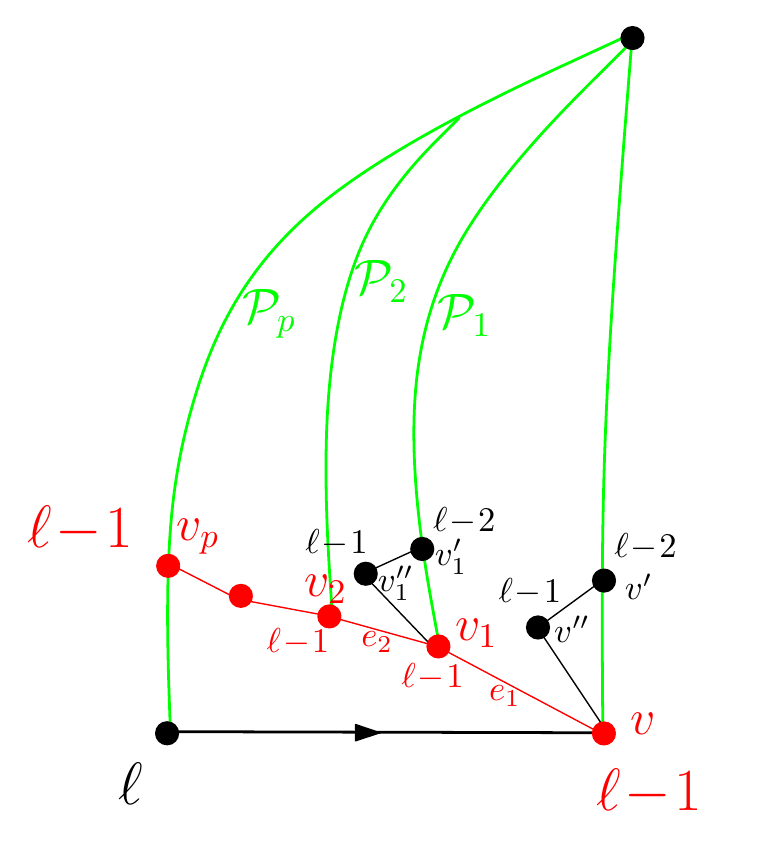}
\end{center}
\caption{Construction of the dividing line (in red) in an $R$-slice (see text).}
\label{fig:dividing}
\end{figure}
\begin{figure}
\begin{center}
\includegraphics[width=11cm]{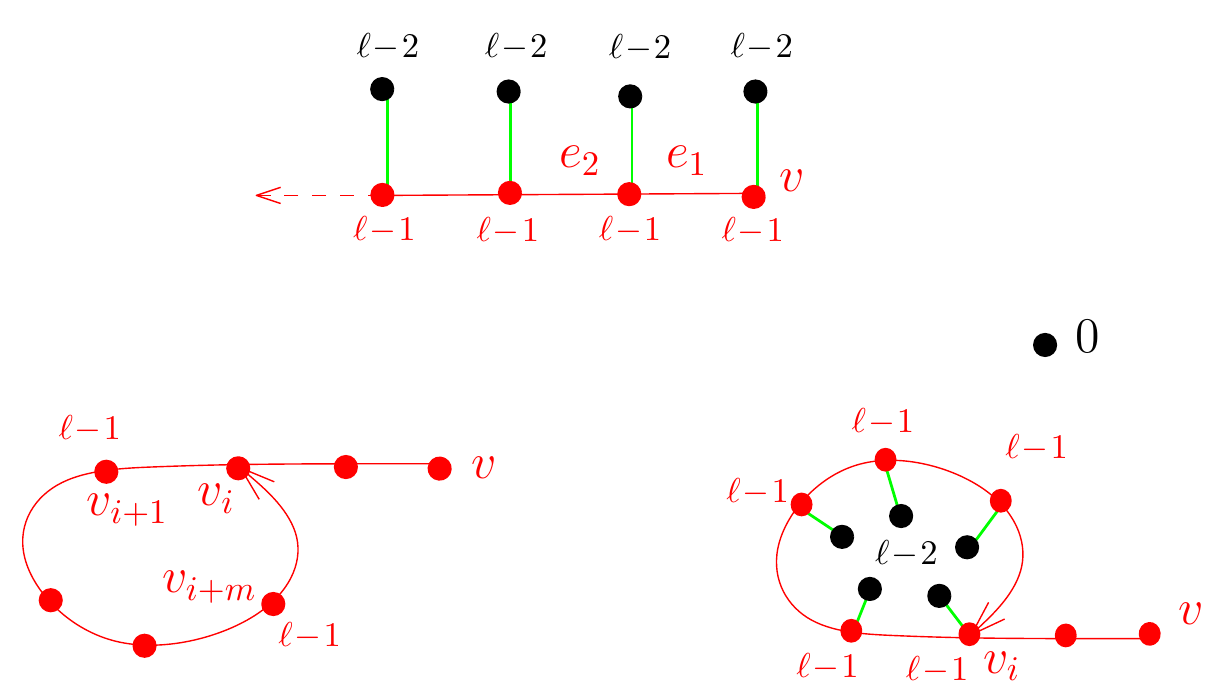}
\end{center}
\caption{Top: a schematic picture of the dividing line. Bottom: if this line forms a loop, we end up with a contradiction at the first revisited vertex $v_i$: 
if the loop closes from the left (bottom left), it means that, at step $i$, we should have picked the edge leading to $v_{i+m}$, not to $v_{i+1}$;
if the loop closes from the right (bottom right), it encloses vertices at distance $\ell-2$ from the apex (marked $0$), which is inconsistent 
with the fact that all the vertices on the dividing line are at distance $\ell-1$ from the apex.}
\label{fig:noloop}
\end{figure}

We shall now derive a new set of recursion relations for $R_k$ and $S_k$ (or more precisely $T_k$) based on a new decomposition of the slices.
This decomposition makes use of a particular \emph{dividing line} drawn on the slice, which we define now. We start for convenience with an $R$-slice whose left boundary has some length $\ell\geq 2$. Its dividing line will then be made of a sequence of edges linking the right and left boundaries of the slice and connecting only
vertices at distance $\ell-1$ from the apex. It is defined recursively as follows (see figure~\ref{fig:dividing}): consider the face on the left of the
first edge of the right boundary, i.e.\ the edge linking the endpoint $v$ of the root-edge, at distance $\ell-1$ from the apex
to its neighbor $v'$ along the right boundary, at distance $\ell-2$ from the apex. The third vertex incident to this face, $v''$, is necessarily different from
$v$ and $v'$ as otherwise, we would have a second edge linking $v$ and $v'$ within the map, hence a second shortest
path from $v$ to the apex, lying strictly on the left of the right boundary. Moreover, $v'$ is necessarily at distance $\ell-1$ from 
the apex: indeed the only allowed values for the distance are $\ell-2$ and $\ell-1$ but the value $\ell-2$ is forbidden as again it would lead
to the existence of another shortest path from $v$ to the apex, strictly on the left of the right boundary. We conclude that $v$ is
incident to 
at least one edge \emph{leading to a distinct neighbor at distance $\ell-1$} from the apex. Let us pick the leftmost such edge $e_1$
and call $v_1$ its endpoint. Assuming that $v_1$ does not belong to the left boundary, we draw the \emph{leftmost} shortest path 
$\mathcal{P}_1$ from $v_1$ to the apex and call
$v_1'$ the vertex on $\mathcal{P}_1$ at distance $\ell-2$ from the apex. Consider again the face on the left of the
first edge of $\mathcal{P}_1$, linking  $v_1$ to $v_1'$. It is incident to a third vertex $v_1''$, necessarily different from
$v_1$ and $v_1'$ and at distance $\ell-1$ from the apex (for the same reasons as above). We thus conclude that $v_1$ is
incident to 
at least one edge leading to a distinct neighbor at distance $\ell-1$ from the apex. As before, we pick the leftmost such edge $e_2$
and call $v_2$ its endpoint. We may repeat the procedure as long as we do not reach the left boundary, thus creating an oriented 
line $(e_1,e_2,\cdots)$ linking only vertices at distance $\ell-1$ from the apex with moreover, on the right of each vertex along the line,
an edge linking this vertex to a vertex at distance $\ell-2$ from the apex (see figure~\ref{fig:noloop}-top).
It is easy to see that this line cannot make a loop.
Indeed, let us assume that the line revisits some already visited vertex and pick the first such vertex. If this vertex is reached from the left, this contradicts the fact
that, in our construction, we always picked the leftmost edge to a neighbor at distance $\ell-1$ from the apex (see figure~\ref{fig:noloop}-bottom left). 
If it is reached from the right,
this creates a closed region surrounded by vertices at distance $\ell-1$ from the apex, which does dot contain the apex and which 
contains vertices at distance $\ell-2$ from the apex, a contradiction (see figure~\ref{fig:noloop}-bottom right). The line thus necessarily ends after a finite number $p$ of
steps at the vertex $v_p$ lying on the left boundary at distance $\ell-1$ from the apex (note that $\mathcal{P}_p$ 
is then the part of the left boundary lying between $v_p$ and the apex). The open line $(e_1,e_2,\cdots,e_p)$ from $v$ to $v_p$ constitutes our
dividing line with the following property: it is a simple curve linking the right and left boundaries, visiting only vertices 
at distance $\ell-1$ from the apex, dividing {\it de facto} the map in two parts, an upper part containing the apex and a lower part containing the 
root-vertex. Finally, by construction, we have the following property:
\begin{proper}{}
\emph{Two vertices of the dividing line cannot be linked by an edge lying strictly inside the lower part. }
\end{proper}

Indeed, violating Property 1 would contradict the fact that, in our construction, we always picked the leftmost edge leading to a neighbor at distance $\ell-1$.
We could similarly have started with an $S$-slice whose left boundary has some length $\ell\geq 2$. The dividing line would then be defined exactly in the same way, now starting from the vertex $v$ of the right boundary at distance $\ell-1$ from the apex (see figure~\ref{fig:Tknewrec} 
for an illustration). 

In the following, we shall recourse to the dividing line to decompose $R$-slices enumerated by $R_k$ and $S$-slices enumerated 
by $T_k$. Both families of slices have a left-boundary length $\ell$ satisfying $1\leq \ell \leq k$. So far
we defined the dividing line only for $\ell \geq 2$. For $\ell=1$, we take the convention that the dividing line is reduced to a single vertex equal to the apex.
\vskip .2cm
\begin{figure}
\begin{center}
\includegraphics[width=7cm]{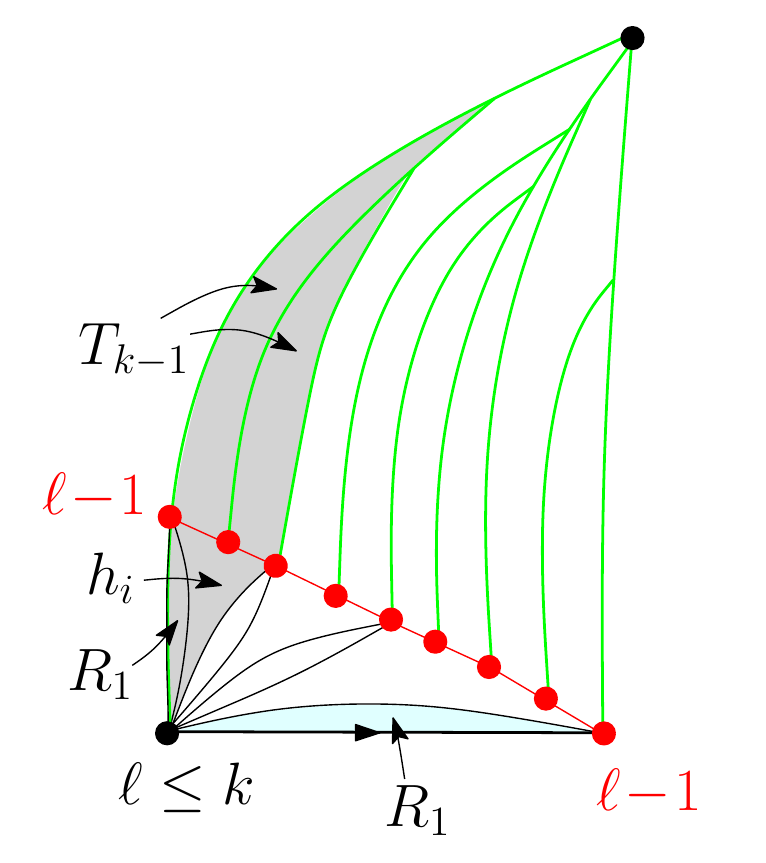}
\end{center}
\caption{The decomposition of an $R$-slice enumerated by $R_k$ into a sequence of blocks. The first block is indicated in gray. Each block
is formed of a bundle (enumerated by $R_1$), a rooted triangulation with a boundary of some arbitrary length 
$i\geq 3$ and with particular properties (see text), enumerated by $h_i$, and a set of $i-2$ attached $S$-slices enumerated by $T_{k-1}$. 
The generating function of a block is thus $R_1\, \sum_{i\geq 3} h_i\, T_{k-1}^{i-2}$. The sequence 
of blocks is to be completed by a final bundle (enumerated by $R_1$, here in light blue) connecting the extremities of the root-edge. }
\label{fig:Rknewrec}
\end{figure}
\begin{figure}
\begin{center}
\includegraphics[width=7cm]{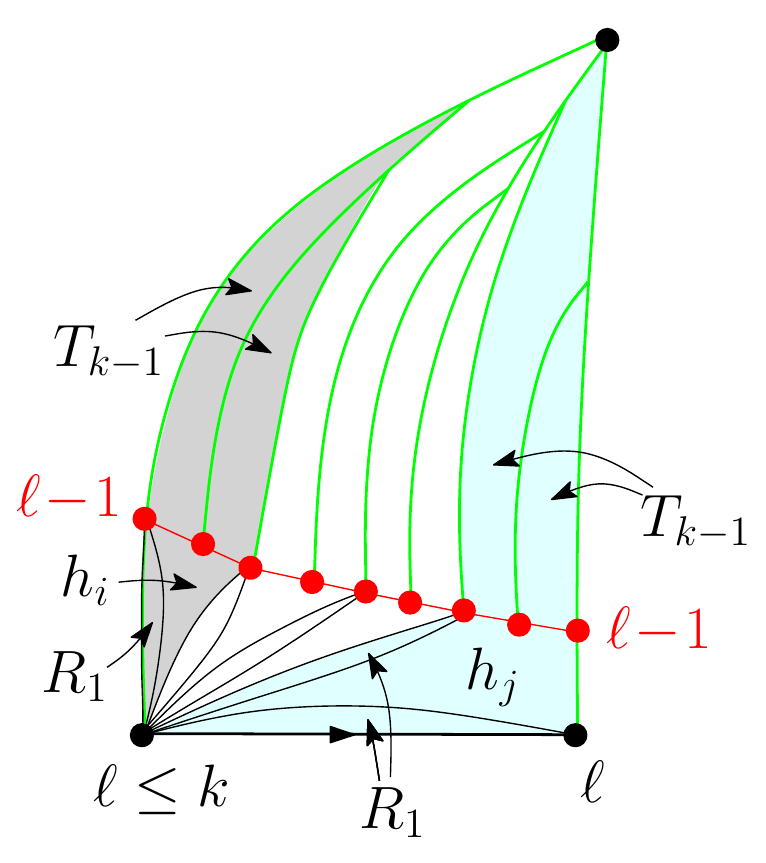}
\end{center}
\caption{The decomposition of an $S$-slice enumerated by $T_k$ into blocks similar to those of figure~\ref{fig:Rknewrec},
each enumerated by $R_1\, \sum_{i\geq 3} h_i\, T_{k-1}^{i-2}$.
The sequence of blocks is now to be completed by a final portion (in light blue) formed of two bundles (each enumerated by $R_1$), a rooted triangulation
with a boundary of some arbitrary length $j\geq 3$, enumerated by $h_j$, and a set of $j-3$ $S$-slices enumerated by $T_{k-1}$.
The generating function of the last portion is thus $R_1^2\, \sum_{j\geq 3} h_j\, T_{k-1}^{j-3}$.}
\label{fig:Tknewrec}
\end{figure}
\subsection{A new set of recursion relations} We shall now derive a new set of recursion relations for $R_k$ and $T_k$,
based on a new decomposition of slices intimately linked to their dividing line. We start again for convenience with an $R$-slice and draw its dividing line, as defined above. We then note that the root-vertex, 
at distance $\ell$ from the apex, is adjacent, in all generality, to a number of vertices of the dividing line. These include the two extremities of the line, plus possibly some of its internal vertices. 
In general, such adjacency with a given vertex of the dividing line is moreover achieved by a bundle of edges, as we defined it, 
which we view as a rooted triangulation with
a boundary of length two (the boundary being made of the two outermost edges performing the connection), hence which is enumerated by $R_1$. For each vertex
of the dividing line adjacent to the root-vertex, we cut the map along the \emph{leftmost} edge of the associated bundle and along the leftmost shortest path
from this vertex to the apex. This cutting decomposes the slice into a sequence of blocks (see figure~\ref{fig:Rknewrec}), each of which being 
formed of (1) a bundle enumerated by $R_1$, 
(2) a triangulation with a boundary of some arbitrary length $i\geq 3$, lying in the lower part of the slice in-between two successive bundles, and enumerated by a generating function $h_i\equiv h_i(g)$ that we shall analyze just below, 
and (3) a set of $i-2$ $S$-slices \emph{whose
root-edge does not form a loop}. If we start with an $R$-slice enumerated by $R_k$, hence with $1\leq \ell \leq k$, then the $\ell=1$ contribution 
yields the empty sequence of blocks (since the dividing line is reduced to the apex in this case), while the $\ell \geq 2$ contribution yields non-empty sequences of blocks whose 
$S$-slice components have some arbitrary left-boundary length between $1$ and $\ell-1$, hence between $1$ and $k-1$, 
as enumerated by $T_{k-1}$ (since their root-edge cannot form a loop). 
To summarize, each block of the (possibly empty) sequence is enumerated by 
$R_1\, \sum_{i\geq 3} h_i\, T_{k-1}^{i-2}$. Finally, we are left with a final bundle connecting the two extremities of the root-edge and enumerated
by $R_1$ (see figure~\ref{fig:Rknewrec}). If we now start instead with an $S$-slice enumerated by $T_k$, hence satisfying $1\leq \ell \leq k$, a similar decomposition produces a (possibly empty) sequence of the same blocks, now completed by a final portion of map enumerated by
$R_1^2\, \sum_{i\geq 3} h_i\, T_{k-1}^{i-3}$ (see figure~\ref{fig:Tknewrec}).
We arrive at the relations
\begin{equation}
\left\{
\begin{split}
R_k& =\frac{R_1}{1-R_1\sum\limits_{i\geq 3}h_i\, T_{k-1}^{i-2}}=\frac{R_1}{1-R_1\, T_{k-1}\, \Phi(T_{k-1})}\\
T_k& =\frac{R_1^2 \sum\limits_{i\geq 3}h_i, T_{k-1}^{i-3}}{1-R_1\sum\limits_{i\geq 3}h_i, T_{k-1}^{i-2}} 
=\frac{R_1^2\, \Phi(T_{k-1})}{1-R_1\, T_{k-1}\, \Phi(T_{k-1})}\\
\end{split}
\right.
\label{eq:newrec}
\end{equation}
with
\begin{equation}
\Phi(T)\equiv \Phi(T,g)=\sum\limits_{i\geq 3}h_i(g)\, T^{i-3}\ .
\label{eq:defPhi}
\end{equation}
Note in particular the relation
\begin{equation}
R_k-R_1=T_{k-1}\, T_k\ .
\label{eq:relRkTk}
\end{equation}
So far we did not discuss the precise definition of $h_i\equiv h_i(g)$ ($i\geq 3$). By construction, the maps enumerated by $h_i$ correspond to 
a part of the slice lying below the dividing line and in-between two consecutive bundles. This part forms a rooted triangulation
with a boundary of length $i$ made of a segment formed by $i-2$ consecutive edges of the dividing line and of $2$ extra edges connecting the extremities of this segment to the root-vertex
(we may then decide for instance to root the map at its rightmost edge from the root vertex to the dividing line, oriented away from the root-vertex).
This boundary forms by construction a simple curve. Moreover, we have the property:
\begin{proper}
\emph{In the maps enumerated by $h_i$, two vertices of the boundary cannot be linked by an edge lying strictly inside the map.}
\end{proper}
For a pair of boundary vertices belonging to the dividing line, this property follows immediately from Property 1 above. As for the root vertex, the only boundary vertices to which 
it is connected are the two extremities of the segment of the dividing line and each of this connection is performed by a single boundary edge.  

To summarize, $h_i\equiv h_i(g)$ is the generating function of rooted triangulations with a boundary of length $i$ ($i\geq 3$) forming a simple curve, 
and \emph{with the property 2} above. It is interesting to note that, even if, in the construction, the root-vertex and the vertices
of the dividing line play very different roles, all boundary vertices eventually play \emph{symmetric roles} in the maps enumerated by $h_i$.
\vskip .2cm
Assuming that we know $\Phi(T)$, the second line of \eqref{eq:newrec} is a direct recursion on $k$ which fixes $T_k$ for all $k\geq 0$
recursively from the initial condition $T_0=0$.
Then the first line of \eqref{eq:newrec} gives access to $R_k$ for all $k\geq 1$. Strictly speaking, we also need as input the knowledge of
$R_1$ and, if we eventually want an expression for $S_k=T_k+S_0$ and for the two-point function $G_k$, 
we also need the value of $S_0$. We will explain in Section \ref{sec:finalexpr} how to get rid of this problem.

\subsection{Back to Tutte's seminal paper}
\label{sec:TutteCPT}
The natural question as this stage is: do we have an expression for $\Phi(T)$? Remarkably, the answer is yes, as shown
in Tutte's seminal paper \cite{TutteCPT} on the enumeration of triangulations. There Tutte introduces precisely the 
same notion of triangulations having a boundary forming a simple curve of arbitrary length at least $3$ and \emph{satisfying Property 2} above.
To be precise, Tutte considers what are called \emph{simple triangulations}, i.e.\ triangulations required to have no loops nor multiple edges.
The quantity considered by Tutte is therefore a slight reformulation of $\Phi(T)$, called $\psi(x,y)$ in \cite{TutteCPT} ,
but the passage from $\psi(x,y)$ to $\Phi(T)$ is straightforward and can be obtained via a simple substitution procedure. This 
correspondence will be made explicit in Section~\ref{sec:simpletriang} below. With this correspondence, Tutte's result immediately translates into
the following equation
\begin{equation}
R_1\, T^2 \, \Phi^2(T)+(g\, R_1^2+g\, R_1^3\, h_3\, T-T-g\, R_1\, T^2)\, \Phi(T)+(g\,T-g\, R_1^2\, h_3)=0
\label{eq:eqPhi}
\end{equation}
which, as shown in \cite{TutteCPT},  \emph{entirely fixes $\Phi(T)$} as a function of $T$, $g$ and $R_1$ (see Section~\ref{sec:Tuttesol} below
for an explicit expression).
We thus have at our disposal all the ingredients to solve our new recursion relation, a task that will be performed
explicitly in the next sections.
\vskip .2cm
\begin{figure}
\begin{center}
\includegraphics[width=6cm]{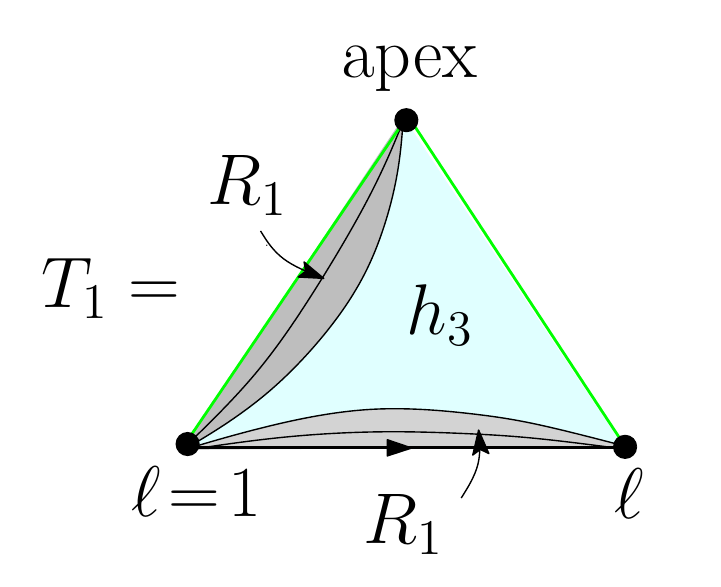}
\end{center}
\caption{A schematic explanation of the relation $T_1=R_1^2\, h_3$.}
\label{fig:T1h3}
\end{figure}

Before we solve our recursion relation, it is interesting to explore what we may learn by simply \emph{making the system \eqref{eq:newrec} consistent
with the more classical system \eqref{eq:oldrec}}. 
The first line in \eqref{eq:newrec} may be rewritten as
\begin{equation*}
T_{k-1}\Phi(T_{k-1})=\frac{1}{R_1}-\frac{1}{R_k}
\end{equation*}
while the first line of \eqref{eq:oldrec} leads to
\begin{equation*}
\frac{1}{R_k}=1-g\, (S_k+S_{k-1})=1-2g\, S_0-g\, (T_k+T_{k-1})
\end{equation*}
and in particular
\begin{equation*}
\frac{1}{R_1}=1-2g\, S_0-g\, T_1=1-2g\, S_0-g\, R_1^2\, h_3
\end{equation*}
since $T_0=0$ and $T_1=h_3 R_1^2$. This later relation is easily understood by noting that $T_1$ enumerates $S$-slices 
with $\ell=1$, which are rooted triangulations with a boundary of length $3$ forming a simple curve (see figure~\ref{fig:T1h3}). 
On the other hand, $h_3$ also enumerates rooted triangulations with a boundary of length $3$ and the only difference 
is that, in $T_1$, both the root edge and the left boundary edge may be doubled by a bundle of edges (the right boundary edge 
cannot be doubled by a bundle as it is the unique shortest path to the apex). In other words, to get $T_1$, we must multiply $h_3$ 
twice by the bundle generating function $R_1$, hence the relation.

Combining the above equations, we deduce
\begin{equation*}
T_{k-1}\Phi(T_{k-1})=g\, (T_k+T_{k-1})-g\, R_1^2\, h_3=g\, \left(\frac{R_1^2 \, \Phi(T_{k-1})}{1-R_1\, T_{k-1}\, \Phi(T_{h-1})}+T_{k-1}\right)-g\, R_1^2
\, h_3
\end{equation*}
where we have used the second line of \eqref{eq:newrec} to express $T_k$ in terms of $T_{k-1}$.
Equating the left and right terms above the  leads precisely to equation \eqref{eq:eqPhi} for the specific value $T=T_{k-1}$. 
This equation being valid for any positive integer $k$, we may reasonably infer that it holds for any $T$ (small enough
so that $\Phi(T)$ is well-defined). Indeed the explicit dependence of $T_{k-1}$ in $k$ (see below) allows to formally extend $T_{k-1}$
to non-integer values of $k$, so that $T_{k-1}$ now varies continuously with real $k$.
To summarize, making the systems \eqref{eq:newrec} and \eqref{eq:oldrec} consistent is yet another way to understand
Tutte's equation \eqref{eq:eqPhi}. 

\section{A detour via simple triangulations}
\label{sec:simpletriang}
\subsection{Substitution}
\label{sec:subst}
As already mentioned, the analysis of \cite{TutteCPT} deals with simple triangulations, i.e.\ triangulations 
\emph{with neither loops nor multiple edges} (here a loop
stands of an edges with identical extremities).
Let us therefore introduce the generating functions $r_k$, $t_k$ ($k\geq 1$) and $\tilde{h}_i$ ($i\geq 3$), defined as $R_k$, $T_k$ and $h_i$
respectively, but with the constraint that the map contains neither loops nor multiple edges. It is easy to see that,
for maps enumerated by $R_k$, $T_k$ and $h_i$, the presence of a loop automatically implies the presence of a multiple edge. 
Indeed, a loop separates the map into two regions, its exterior, which contains the outer face and its interior. One loop is then said to be included 
in another if it lies in its interior. This inclusion defines a partial ordering of the loops and we may consider one of the largest elements for this ordering,
i.e.\ a loop which is not contained in the interior of any other loop. The face incident to the edge forming this loop and lying
in its exterior is necessarily a triangle of the bulk. Indeed, the boundary of maps enumerated by $R_k$, $T_k$ or $h_i$ cannot contain loops.
The two remaining edges of this triangle necessarily form a multiple edge surrounding the loop by connecting the endpoint of the loop to a distinct
vertex in the exterior of the loop (the other possibility, namely that the two remaining edges of the triangle form two loops, is ruled out
as, if so, one of these two loops would encircle the supposedly largest loop, a contradiction).

To suppress both loops and multiple edges in maps enumerated by $R_k$, $T_k$ and $h_i$, it is thus sufficient to 
\emph{suppress multiple edges only}.
This translates into a simple substitution procedure to go from $r_k$, $t_k$, $\tilde{h}_i$ to $R_k$, $T_k$ and $h_i$:
maps in the second family are directly obtained by simply replacing the edges in the maps of the first family by bundles 
of edges, as enumerated by $R_1$. Note that the substitution is to be performed for each edge except for some of the
boundary edges: in the case of $r_k$ and $t_k$, edges of the right boundary must be left untouched as duplicating
some of them would create a shortest path strictly on the left of the right-boundary. In the case of $\tilde{h}_i$,
none of the boundary edges can be duplicated since we have to enforce Property 2 on maps enumerated by $h_i$.

Let us now describe in details the consequences of the substitution procedure. Consider a simple triangulation with a simple boundary of length 
$L$ and call $E$, $V$ and $F$ its numbers of \emph{inner} edges,
vertices and \emph{inner} faces respectively.
From the relation $3F=2E+L$ and the Euler relation $F+V-E-L=1$, we immediately deduce
\begin{equation}
V=\frac{L+F}{2}+1\ ,\qquad E=\frac{3F-L}{2}\ .
\label{eq:Euler}
\end{equation}

\noindent{\underline{\emph {The case of $h_i$ vs $\tilde{h}_i$}}}\,:\ In this case, we have $L=i$ and the substitution requires a weight $R_1$ per 
inner edge of the simple triangulation. These edges are in number $(3F-i)/2=3F/2-i/2$, hence we should give an extra weight
$R_1^{3/2}$ per face, resulting in a total face weight
\begin{equation}
G=g\, R_1^{3/2}
\label{eq:valG}
\end{equation} 
together with a global factor $R_1^{-i/2}$. In other words, we have
\begin{equation}
h_i (g)=R_1^{-i/2}\, {\tilde h}_i(G)\ .
\label{eq:valhell}
\end{equation}

\noindent{\underline{\emph{The case of $T_k$ vs $t_k$}}}\,:\ For a slice enumerated by $t_k$, of arbitrary left-boundary length $\ell$ (between $1$ and $k$),  
we have $L=2\ell+1$ and the substitution requires a weight $R_1$ for each inner edge of the simple slice, as well
as for the base edge and for the edges of the left boundary (as already mentioned, there is no substitution attached
to the edges of the right boundary as this boundary must be the unique shortest path between its extremal vertices). 
These edges are in number $(3F-(2\ell+1))/2+\ell+1=3F/2+1/2$, hence we should give an extra weight
$R_1^{3/2}$ per face as before, hence a total weight $G$,  
together with a global factor $R_1^{1/2}$. In other words, we now have
\begin{equation*}
T_k (g)=R_1^{1/2}\, t_k(G)\ .
\end{equation*}

\noindent{\underline{\emph{The case of $R_k$ vs $r_k$}}}\,:\ For a slice enumerated by $r_k$, of arbitrary left-boundary length $\ell$ (between $1$ and $k$),  
we have $L=2\ell$ and the substitution requires a weight $R_1$ for each inner edge of the simple slice, as well
as for the base edge and for the edges of the left boundary. 
These edges are in number $(3F-2\ell)/2+\ell+1=3F/2+1$, hence we should again give a total weight $G$
to each face,  
together with a global factor $R_1$. In other words, we have
\begin{equation*}
R_k (g)=R_1\, r_k(G)\ .
\end{equation*}
\vskip .2cm
Introducing
\begin{equation*}
\tilde{\Phi}(t)\equiv \tilde{\Phi}(t,G)=\sum_{i\geq 3} \tilde{h}_i(G)\, t^{i-3}\ ,
\end{equation*}
we read from \eqref{eq:valhell} the correspondence
\begin{equation}
\Phi (T)=R_1^{-3/2} \tilde{\Phi}(t)\ ,\qquad T=R_1^{1/2}\, t\ .
\label{eq:PhiPhitilde}
\end{equation}
With the above correspondence, the system \eqref{eq:newrec} is equivalent to the relations
\begin{equation}
\left\{
\begin{split}
r_k& =\frac{1}{1-\sum\limits_{i\geq 3} \tilde{h}_i\, t_{k-1}^{i-2}}=\frac{1}{1-t_{k-1}\, \tilde{\Phi}(t_{k-1})}\\
t_k& =\frac{\sum\limits_{i\geq 3}\tilde{h}_i\, t_{k-1}^{i-3}}{1-\sum\limits_{i\geq 3}\tilde{h}_i\, t_{k-1}^{i-2}} 
=\frac{\tilde{\Phi}(t_{k-1})}{1-t_{k-1}\, \tilde{\Phi}(t_{k-1})}\\
\end{split}
\right.
\label{eq:newrecsimple}
\end{equation}
which determine $r_k$ and $t_k$ recursively from the initial condition $t_0=0$, while \eqref{eq:relRkTk} becomes
\begin{equation}
r_k=1+t_{k-1}\, t_k\ .
\label{eq:relrktk}
\end{equation}
Note that all these latter equations could have been obtained directly by applying the decompositions
used in Section~\ref{sec:newrec} directly to maps enumerated by $r_k$ and $t_k$. Note also that 
considering simple triangulations is a way to get rid of $R_1$ (as well as of $S_0$),
which disappeared from our recursion relations. As a final remark,
it would be tempting to believe that the distance-dependent two-point function of simple triangulations is
given by a formula as simple as eq.~\eqref{eq:twopoint}, say by simply replacing $R_k$ and $S_k$ by $r_k$ and $t_k$.
This is however not true since, in the cutting procedure illustrated in figures~\ref{fig:twopoint2} and \ref{fig:twopoint1}, the requirement 
of having no loop nor multiple edge \emph{encircling the marked vertex} introduces non-trivial constraints on the 
associated slices, which are difficult to handle. So our detour in the ensemble of simple triangulations should here be viewed as a simple 
trick to simplify our recursion and to directly use the results of \cite{TutteCPT}. 

\subsection{Equation for $\mathbf{\tilde{\Phi}(t)}$}
\label{sec:eqforPhitilde}
With the correspondence \eqref{eq:valG}, \eqref{eq:valhell} and \eqref{eq:PhiPhitilde}, eq.~\eqref{eq:eqPhi} is fully equivalent to the following equation for $\tilde{\Phi}$:
\begin{equation}
t^2 \, \tilde{\Phi}^2(t)+(G+G\, \tilde{h}_3\, t-t-G\, t^2)\, \tilde{\Phi}(t)+(G\,t-G\, \tilde{h}_3)=0
\label{eq:eqPhitilde}
\end{equation}
where $\tilde{h}_3=R_1^{3/2}\, h_3$ from \eqref{eq:valhell}.
\begin{figure}
\begin{center}
\includegraphics[width=6cm]{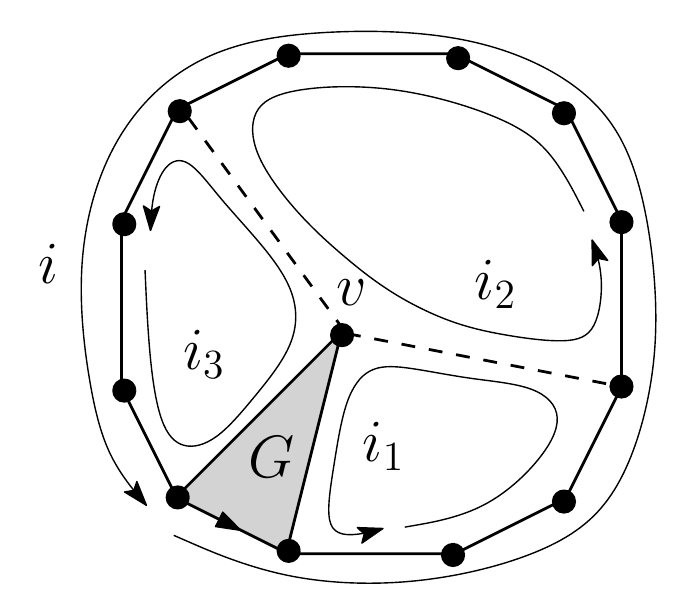}
\end{center}
\caption{The decomposition of a map enumerated by $h_i$ ($i\geq 3$) into $(p+1)$ domains (here $p=2$) enumerated by $h_{i_1}$, $h_{i_2}$, $\cdots ,h_{i_{p+1}}$
respectively. The boundary lengths of the domains satisfy $i_m\geq 3$ for $m=1,\cdots ,p+1$ and $\sum_{m=1}^{p+1}(i_m-2)=i-1$.}
\label{fig:Tuttedecomp}
\end{figure}
Let us recall here how to derive this equation, following Tutte's argument in \cite{TutteCPT}.
Consider a map enumerated by $\tilde h_i$ ($i\geq 3$): for $i=3$, the bulk of the map may possibly be reduced to a single triangle, 
contributing $G$ to $\tilde{h}_3$. In all the other cases, in order to guarantee Property 2, the triangle on the left of the root-edge of the map is incident to 
a vertex $v$ lying strictly inside the bulk (see figure~\ref{fig:Tuttedecomp}). This vertex is in all generality connected to a number $p\geq 0$ of
the $i-2$ boundary vertices different from the extremities of the root-edge. These connections are moreover performed by $p$ single edges.
Removing the triangle on the left of the root-edge and cutting along these $p$ single edges gives rise to $p+1$ domains which are triangulations
with boundaries of lengths $i_1,i_2,\cdots,i_{p+1}$, all larger than $3$, satisfying 
\begin{equation*}
\sum_{m=1}^{p+1}(i_m-2)=i-1\ .
\end{equation*}
The boundaries of these triangulations are simple curves and two boundary vertices cannot be linked by an internal edge. So the 
$m$-th piece is enumerated by $\tilde{h}_{i_m}$. This leads to the identity
\begin{equation*}
\tilde{h}_i=G\, \delta_{i,3}+G\, \sum_{p\geq 0}\sum_{{(i_1,\cdots,i_{p+1})\atop i_m\geq 3\ , m=1,\cdots, p+1}\atop
\sum\limits_{m=1}^{p+1}(i_m-2)=i-1} \tilde{h}_{i_1}\cdots \tilde{h}_{i_{p+1}}\ .
\end{equation*}
Multiplying by $t^{i-3}$ and summing over $i\geq 3$, this rewrites as 
\begin{equation*}
\begin{split}
\tilde{\Phi}(t)&=G+\frac{G}{t^2}\, \sum_{p\geq 0}\sum_{{(i_1,\cdots,i_{p+1})\atop i_m\geq 3\ , m=1,\cdots, p+1}\atop
\sum\limits_{m=1}^{p+1}(i_m-2)\geq 2} (\tilde{h}_{i_1}\, t^{i_1-2})\cdots (\tilde{h}_{i_{p+1}}\, t^{i_{p+1}-2})\\
&= G+\frac{G}{t^2}\, \left(\sum_{p\geq 0}\left(\sum_{i_\geq 3} \tilde{h}_i\, t^{i-2}\right)^{p+1}-\tilde{h}_3\, t\right)\ .\\
\end{split}
\end{equation*}
Note the subtracted term corresponding to $p=0$ and $i_1=3$, which, for arbitrary $p\geq 0$ and $i_m\geq 3$ is
the only case for which the condition $\sum_{m=1}^{p+1}(i_m-2)\geq 2$ is not satisfied.
We end up with 
\begin{equation*}
\begin{split}
\tilde{\Phi}(t) &=  G+\frac{G}{t^2}\, \left(\frac{t\, \tilde{\Phi}(t)}{1-t\, \tilde{\Phi}(t)}-\tilde{h}_3\, t\right)
=  G+\frac{G}{t}\, \left(\frac{\tilde{\Phi}(t)}{1-t\, \tilde{\Phi}(t)}-\tilde{h}_3\right) \\
\end{split}
\end{equation*}
which immediately leads to \eqref{eq:eqPhitilde}, and after substitution to the announced equation \eqref{eq:eqPhi}.

\section{Using Tutte's solution}
\label{sec:Tuttesol}
\subsection{Tutte's generating function $\mathbf{\psi(x,y)}$}
\label{sec:Tutteeq}
We shall now rely on Tutte's solution of the equation \eqref{eq:eqPhitilde} to solve our new recursion relation.
In order to directly use Tutte's expressions in \cite{TutteCPT}, we need a slight (and harmless) reformulation of the
generating function $\tilde{\Phi}(T)$. First, as noted in \cite{TutteCPT} for triangulations enumerated
by $\tilde{h}_i$, the numbers $E$, $V$, $F$ and $L=i$ of, respectively, inner edges, vertices,  inner faces and boundary-edges 
(satisfying \eqref{eq:Euler}) may be written as
\begin{equation*}
E=3n+m\ , \qquad
V=n+m+3\ , \qquad
F=2n+m+1\ , \qquad
L=i=m+3
\end{equation*}
for some $m,n\geq 0$ (since $i\geq 3$ and $V\geq i$).
Using the variables $m$ and $n$ (instead of $F$ and $i$), Tutte introduces (instead of $\tilde{\Phi}(t)=\tilde{\Phi}(t,G)$) the generating function
\begin{equation*}
\psi(x,y)\equiv \sum_{m,n\geq 0}\psi_{m,n}\, x^n\, y^m\ , \qquad \psi_{m,n}\equiv [G^{2n+m+1}]\tilde{h}_{m+3}\ .
\end{equation*}
We immediately read the correspondence
\begin{equation*}
\begin{split}
\psi(x,y)&=\sum_{m\geq 0} y^m \sum_{n\geq 0} x^n [G^{2n+m+1}]\tilde{h}_{m+3}\\
&= \sum_{m\geq 0} y^m\,  \frac{\tilde{h}_{m+3}(G)}{G^{m+1}} \qquad  \quad \hbox{for}\ x=G^2 \\
&= \frac{1}{G} \sum_{m\geq 0} t^m\, \tilde{h}_{m+3}(G)\qquad  \hbox{for}\ y=G\, t \\
&=\frac{\tilde{\Phi}(t)}{G}\\
\end{split}
\end{equation*}
or, in short
\begin{equation*}
\psi(x,y)=\frac{\tilde{\Phi}(t)}{G}\ , \qquad x=G^2\ , \qquad y=G\, t\ .
\end{equation*}
Note in particular that, setting $t=0$, we have 
\begin{equation*}
g_3(x)\equiv \psi(x,0)=\frac{\tilde{h}_3}{G}\ .
\end{equation*}
Setting the correspondence 
\begin{equation}
\begin{gathered}
\underbrace{
y=G\, t\ , \qquad x= G^2\, \qquad g_3=\tilde{h}_3/G\, \qquad \psi=\tilde{\Phi}/G }\\
\Updownarrow \\
\overbrace{t=y/\sqrt{x} \ , \qquad  G=\sqrt{x}\ , \qquad \tilde{h}_3=g_3\, \sqrt{x} \ , \qquad \tilde{\Phi}=\psi\, \sqrt{x} }\\
\end{gathered}
\label{eq:correspondence} 
\end{equation}
equation \eqref{eq:eqPhitilde} is equivalent to
\begin{equation}
y^2\, \psi^2(x,y)+(x+x\, y\, g_3(x)-y-y^2)\psi(x,y)+y-x\, g_3(x)\ , \qquad g_3(x)=\psi(x,0)
\label{eq:eqforpsi}
\end{equation}
which is precisely the form given by Tutte in \cite{TutteCPT}. 

As explained in  \cite{TutteCPT}, the solution of this equation is best expressed upon parametrizing $x$ and $y$ as
\begin{equation*}
\left\{
\begin{matrix}
x=\theta\, (1-\theta)^3 \\  y=(1-\theta)^3\,  \sigma \ .\\
\end{matrix}
\right.
\end{equation*}
With this parametrization, we have \cite{TutteCPT}
\begin{equation*}
x\, g_3(x)= \theta\, (1-2\theta)
\end{equation*}
while $\psi(x,y)$ is fixed by
\begin{equation}
\psi(x,y)=\frac{1}{(1-\theta)^3\,  \sigma}\left(\frac{\theta\, \sigma}{Y(\theta,\sigma)(1+Y(\theta,\sigma))^2}+1\right)
\label{eq:psiformula}
\end{equation}
where $Y(\theta,\sigma)$ is the solution of the quadratic equation
\begin{equation}
Y^2(\theta,\sigma)+ (1-\sigma+\theta\, \sigma)\, Y(\theta,\sigma)+\theta\, \sigma=0
\label{eq:defY}
\end{equation} 
such that $Y(\theta,\sigma)\sim -\theta\, \sigma$ for small $\theta$ or small $\sigma$.

\subsection{Writing the recursion in terms of Tutte's variable}
\label{sec:Tuttevar}
Using the correspondence \eqref{eq:correspondence}, $\theta$ and $\sigma$ are to be considered as parametrizations of $G$ and $t$ 
via
\begin{equation*}
\left\{
\begin{matrix}
G^2=\theta\, (1-\theta)^3 \\  G\, t=(1-\theta)^3\,  \sigma \
 .\\ 
\end{matrix}
\right.
\end{equation*}
In particular, eqs.~\eqref{eq:psiformula} and \eqref{eq:defY} translate into 
\begin{equation*}
\left\{
\begin{split}
\tilde{\Phi}(t)&=\frac{C^3}{Y(t)(1+Y(t))^2}+\frac{1}{t}\\
Y(t)&=\frac{1}{2} \left( C\, t-1+\sqrt{(C\, t-1)^2-4\, C^3\, t} \right)\\
C&=\sqrt{\frac{\theta}{1-\theta}}=\sum_{n\geq 0}\frac{2}{3n+2} {4n+1\choose n}\, G^{2n+1}\ ,\\
\end{split}
\right.
\end{equation*} 
from which (together with the relation $G=g\, R_1^{3/2}$) we may obtain, as announced, an explicit expression of $\Phi(T)=R_1^{-3/2}{\tilde \Phi}(R_1^{-1/2}\, T)$ as a function of $g$ and $R_1$ only. 

We give here the explicit expression of $\tilde{\Phi}(t)$ only for completeness.
To solve our recursion, it is indeed much simpler to directly work with Tutte's variables. In order to write the relations \eqref{eq:newrecsimple}, we must specialize $t$ to the values $t_k$ (and $t_{k-1}$) hence we define
\begin{equation*}
\left\{
\begin{split}
\sigma_{k}&\equiv \frac{G\, t_k}{(1-\theta)^3}=\sqrt{\frac{\theta}{(1-\theta)^3}}\ t_k\\
Y_k& \equiv Y(\theta,\sigma_k)\\
\end{split}
\right.
\end{equation*}  
for $k\geq 0$. Note that, by inverting the relation \eqref{eq:defY} defining $Y(\theta,\sigma)$, we may write $\sigma_k$ in terms of
$Y_k$ as
\begin{equation}
\sigma_k=\frac{Y_k\, (1+Y_k)}{(1-\theta)\, Y_k-\theta}\ .
\label{eq:sigmaY}
\end{equation}  
Let us now show that our recursion relation translates into a \emph{particularly simple recursion relation for $Y_k$}. 
Writing the second line of eq.~\eqref{eq:newrecsimple} in terms of $\sigma_k$,
we get immediately
\begin{equation}
\sigma_k  
=\frac{\theta\, \psi_{k-1}}{1-\sigma_{k-1}\, (1-\theta)^3\,  \psi_{k-1}}\ , \qquad \psi_{k-1}\equiv \psi\Big(\theta\, (1-\theta)^3,(1-\theta)^3\, \sigma_{k-1}\Big)\ .
\label{eq:sigmak}
\end{equation}
Using the relation \eqref{eq:sigmaY} and the expression \eqref{eq:psiformula}, we have the expressions
\begin{equation*}
\left\{
\begin{split}
\sigma_{k-1}&=\frac{Y_{k-1}\, (1+Y_{k-1})}{(1-\theta)\, Y_{k-1}-\theta}\\
\psi_{k-1}&=\frac{1}{(1-\theta)^3\,  \sigma_{k-1}}\left(\frac{\theta\, \sigma_{k-1}}{Y_{k-1}(1+Y_{k-1})^2}+1\right)
=\frac{(1-2\theta)+(1-\theta)\, Y_{k-1}}{(1-\theta)^3(1+Y_{k-1})^2}\\
\end{split}
\right.
\end{equation*}
which, incorporated in \eqref{eq:sigmak} lead to the expression of $\sigma_k$ in terms of $Y_{k-1}$
\begin{equation*}
\sigma_k=\frac{(\theta-(1-\theta)\, Y_{k-1})((1-2\theta)+(1-\theta)\, Y_{k-1})}{(1-\theta)^3\, (1+Y_{k-1})}\ .
\end{equation*}
Comparing with \eqref{eq:sigmaY}, this yields the relation
\begin{equation*}
\frac{Y_k\, (1+Y_k)}{(1-\theta)\, Y_k-\theta}=\frac{(\theta-(1-\theta)\, Y_{k-1})((1-2\theta)+(1-\theta)\, Y_{k-1})}{(1-\theta)^3\, (1+Y_{k-1})}
\end{equation*}
which we may equivalently write as
\begin{equation*}
\left(Y_k+Y_{k-1}+\frac{1-2\theta}{1-\theta}\right)\left(Y_k+\frac{\theta}{(1-\theta)^2}\ \frac{\theta-(1-\theta)\, Y_{k-1}}{1+Y_{k-1}}\right)=0\ .
\end{equation*}
In order to choose which of the two factors we should cancel, we recall that both $Y_k$ and $Y_{k-1}$ should vanish for $\theta \to 0$,
in which case only the second factor above vanishes. 
We are thus led to cancel the second factor in the above product, hence  
we eventually end up with the desired recursion relation for $Y_k$:
\begin{equation}
Y_{k}=-\frac{\theta}{(1-\theta)^2}\ \frac{\theta-(1-\theta)\, Y_{k-1}}{1+Y_{k-1}}\ .
\label{eq:Ykrec}
\end{equation}
This relation is equivalent to our initial recursion \eqref{eq:newrec} for $T_k$. It fixes $Y_k$ for all $k\geq 0$ from the initial condition 
$Y_0=0$ (since $t_0=0$, hence $\sigma_0=0$) and the knowledge of $Y_k$ allows us to obtain $\sigma_k$, $t_k$ and eventually $T_k$.
\subsection{Solving the recursion relation}
\label{sec:solrel}
Getting the solution of the recursion relation \eqref{eq:Ykrec} is a standard exercise and goes as follows: consider more generally
the equation 
\begin{equation*}
Y_k=f(Y_{k-1})\, \qquad f(Y)\equiv \frac{a\, Y+b}{c\, Y+d}\ .
\end{equation*}
Introducing the two fixed points $\alpha$ and $\beta$ of the function $f$ (i.e.\ the two solutions of $f(Y)=Y$),
then the quantity 
\begin{equation*}
W_k=\frac{Y_k-\alpha}{Y_k-\beta}
\end{equation*}
is easily seen to satisfy $W_k=\lambda\, W_{k-1}$, hence
\begin{equation*}
W_k=\lambda^k\, W_0\ , \qquad \lambda\equiv \frac{c\, \beta+d}{c\, \alpha+d}\ .
\end{equation*}
This immediately yields $Y_k$ 
via $Y_k=(\alpha-\beta\, W_k)/(1-W_k)$ (strictly speaking, we must have $\alpha\neq \beta$, which can be verified
a posteriori in our case).

To solve eq.~\eqref{eq:Ykrec}, we may take
\begin{equation*}
a=\frac{\theta}{1-\theta}\ , \qquad b=-\frac{\theta^2}{(1-\theta)^2}\ , \qquad c=1\ , \qquad d=1
\end{equation*}
and thus
\begin{equation*}
f(Y)-Y\propto Y^2+\frac{1-2\theta}{1-\theta}Y+\frac{\theta^2}{(1-\theta)^2}=(Y-\alpha)(Y-\beta)
\end{equation*}
so that 
\begin{equation*}
\alpha+\beta= -\frac{1-2\theta}{1-\theta}\ , \qquad \alpha\, \beta=\frac{\theta^2}{(1-\theta)^2}\ .
\end{equation*}
We deduce
\begin{equation*}
\frac{\lambda}{(1+\lambda)^2}=\frac{\frac{\beta+1}{\alpha+1}}{\left(1+\frac{\beta+1}{\alpha+1}\right)^2}
=\frac{1+\alpha+\beta+\alpha\, \beta}{(2+\alpha+\beta)^2}=\frac{1-\frac{1-2\theta}{1-\theta}+\frac{\theta^2}{(1-\theta)^2}}{(2-\frac{1-2\theta}{1-\theta})^2}
=\theta\ .
\end{equation*}
In other word, we have the correspondence between $\theta$ and $\lambda$
\begin{equation*}
\theta=
\frac{\lambda}{(1+\lambda)^2}=
\frac{1}{\left(\frac{1}{\sqrt{\lambda}}+\sqrt{\lambda}\right)^2}
\end{equation*}
(note that since $\lambda=(\beta+1)/(\alpha+1)$, the condition $\alpha\neq \beta$ is equivalent to the condition $\lambda\neq 1$).
In terms of $\lambda$, we have
\begin{equation*}
\alpha+\beta= -\frac{1+\lambda^2}{1+\lambda+\lambda^2}\ , \qquad \alpha\, \beta=\frac{\lambda^2}{(1+\lambda+\lambda^2)^2}
\end{equation*}
so that we may take
\begin{equation*}
\alpha=-\frac{\lambda^2}{1+\lambda+\lambda^2}\ , \qquad
\beta=-\frac{1}{1+\lambda+\lambda^2}\ .
\end{equation*}
Since $Y_0=0$, we have $W_0=\alpha/\beta=\lambda^2$, so that 
\begin{equation*}
W_k=\lambda^{k+2}\ , \qquad Y_k=\frac{\alpha-\beta\, W_k}{1-W_k}=-\frac{\lambda^2}{1+\lambda+\lambda^2}\times \frac{1-\lambda^k}{1-\lambda^{k+2}}\ .
\end{equation*}
\begin{figure}
\begin{center}
\includegraphics[width=6cm]{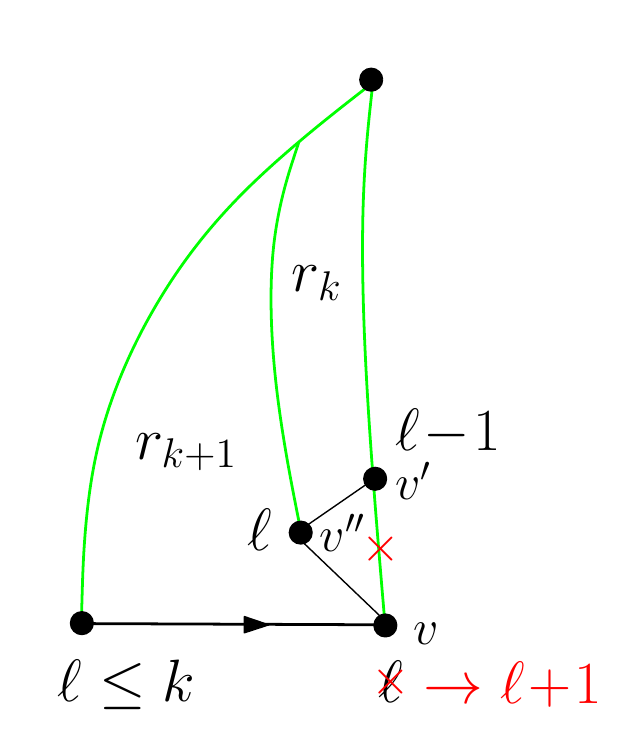}
\end{center}
\caption{A schematic explanation of the relation $t_k=G\, r_k\, r_{k+1}$.}
\label{fig:Grr}
\end{figure}
This solves our recursion relation \eqref{eq:newrecsimple}. Returning to the slice generating functions, we indeed have
\begin{equation*}
\begin{split}
t_k& =\sqrt{\frac{(1-\theta)^3}{\theta}}\, \sigma_k= \sqrt{\frac{(1-\theta)^3}{\theta}}\, \frac{Y_k\, (1+Y_k)}{(1-\theta)\, Y_k-\theta}\\
&= \sqrt{\frac{\lambda}{1+\lambda+\lambda^2}}\times \frac{(1-\lambda^k)(1-\lambda^{k+3})}{(1-\lambda^{k+1})(1-\lambda^{k+2})}\\
\end{split}
\end{equation*}
and from \eqref{eq:relrktk}
\begin{equation*}
r_k=1+t_{k-1}\, t_k=\frac{(1+\lambda)^2}{1+\lambda+\lambda^2}\times \frac{(1-\lambda^k)(1-\lambda^{k+2})}{(1-\lambda^{k+1})^2}
\end{equation*}
where $\lambda$ may be viewed as parametrizing $G$ via
\begin{equation}
G=\sqrt{\theta\, (1-\theta)^3}= \frac{\sqrt{\lambda\, (1+\lambda+\lambda^2)^3}}{(1+\lambda)^4}\ .
\label{eq:Glambda}
\end{equation}
The condition $\lambda \neq 1$ limits the range of $G$ between $0$ and $3\sqrt{3}/16$, in agreement with the fact that the number
of simple triangulations with $n$ faces (and a boundary of finite length) has a large $n$ exponential growth of the form 
$(16/(3\sqrt{3}))^n$ \cite{TutteCPT}.
With their explicit values, we may easily verify the relation
\begin{equation}
t_k=G\, r_k\, r_{k+1}
\label{eq:tGrr}
\end{equation}
which, for $k\geq 1$, can be explained combinatorially as follows:
take a map enumerated by $t_k$, with left-boundary length $\ell$ ($1\leq \ell \leq k$) and consider the triangle immediately on the left of the first
edge of the right boundary linking the endpoint $v$ of the root-edge to a vertex $v'$ at distance $\ell-1$ from the apex
(see figure~\ref{fig:Grr}). The third
vertex $v''$ incident to this triangle is distinct from the other two (as there are no loops) and is at distance $\ell$ from the apex. More
generally, all neighbors of $v$ but $v'$ must be at distance at least $\ell$ from the apex as otherwise, we would have 
a shortest path from $v$ to the apex different from the right-boundary. Removing the edge from $v$ to $v'$ (and the incident triangle), 
the distance from the vertex $v$ to the apex therefore becomes $\ell+1$. Cutting the resulting map along the leftmost shortest path 
from $v''$ to the apex
creates two slices, one enumerated by $r_k$ and the other by $r_{k+1}$ (see figure~\ref{fig:Grr}) hence the relation \eqref{eq:tGrr}.
\section{Final expressions}
\label{sec:finalexpr}
Let us now return to our original problem and obtain expressions for $R_k$, $S_k$ and eventually for the 
two-point function $G_k$.
From the relations 
\begin{equation}
r_k=\frac{R_k}{R_1} \ , \qquad
t_k=\frac{S_k-S_0}{R_1^{1/2}}\ , \qquad
G=g\, R_1^{3/2}
\label{eq:rtG}
\end{equation}
we can get $R_k$ and $S_k$ as functions of $\lambda$, $R_1$ and $S_0$.
Introducing $R_\infty\equiv\lim_{k\to \infty}R_k$, $S_\infty\equiv\lim_{k\to \infty}S_k$
and similarly $r_\infty$ and $t_\infty$, we can write instead $R_k$ and $S_k$ as functions of
$\lambda$, $R_\infty$ and $S_\infty$ via the correspondence
\begin{equation}
\left\{
\begin{split}
r_\infty&=\frac{(1+\lambda)^2}{1+\lambda+\lambda^2}=\frac{R_\infty}{R_1} \\
t_\infty&=\sqrt{\frac{\lambda}{1+\lambda+\lambda^2}}=\frac{S_\infty-S_0}{R_1^{1/2}}\ .\\
\end{split}
\right.
\label{eq:rtinf}
\end{equation}
Eq.~\eqref{eq:rtG} reads indeed
\begin{equation*}
\left\{
\begin{split}
R_k&=R_\infty\, \frac{r_k}{r_\infty} \\
S_k&=S_\infty-\sqrt{\frac{R_\infty}{r_\infty}}\, t_\infty\, \left(1-\frac{t_k}{t_\infty}\right)\\
g&=G\, \left(\frac{r_\infty}{R_\infty}\right)^{3/2}\\
\end{split}
\right.
\end{equation*}
with, as a consequence of \eqref{eq:tGrr}, $t_\infty=G\, r^2_\infty$, hence $\sqrt{R_\infty/r_\infty}\, t_\infty=
g\, R_\infty^2$. We end up with the explicit relations
\begin{equation*}
\left\{
\begin{split}
R_k&=R_\infty\, \frac{(1-\lambda^k)(1-\lambda^{k+2})}{(1-\lambda^{k+1})^2} \\
S_k&=S_\infty-g\, R_\infty^2\, \lambda^k\, \frac{(1-\lambda)(1-\lambda^{2})}{(1-\lambda^{k+1})(1-\lambda^{k+2})}\\
\end{split}
\right.
\end{equation*}
which reproduce the formulas found in \cite{PDFRaman,BG12}. To end our calculation, we still have to express $\lambda$, $R_\infty$ and $S_\infty$
in terms of the weight $g$ only. The quantities $R_\infty$ and $S_\infty$ are simply obtained as the solutions of the 
system obtained by letting $k\to \infty$ in eq.~\eqref{eq:oldrec}, namely
\begin{equation}
\left\{
\begin{split}
R_\infty& =1+2g\, R_\infty\, S_\infty\\
S_\infty& = g\, (S_\infty^2 +2R_\infty)\ .\\
\end{split}
\label{eq:recRS}
\right.
\end{equation}
The desired solution is entirely determined by the condition $R_\infty=1+O(g^2)$ and $S_\infty=2g+O(g^3)$. 
As for $\lambda$, we note that, 
from the expression \eqref{eq:Glambda} of $G$ and that, \eqref{eq:rtinf}, of $r_\infty$, we can immediately
write that 
\begin{equation*}
G^2\, r_\infty^3=\frac{\lambda}{(1+\lambda)^2}=\frac{1}{\lambda+\frac{1}{\lambda}+2}
\end{equation*}
while, from eqs.~\eqref{eq:rtG} and \eqref{eq:rtinf}, $G^2\, r_\infty^3=g^2\, R_\infty^3$. To summarize, $\lambda$ 
is connected to $g$ via
\begin{equation}
\lambda+\frac{1}{\lambda}+2=\frac{1}{g^2\, R_\infty^3}
\label{eq:charac}
\end{equation}
which is precisely the relation found in \cite{BG12}. To be as explicit as in the case of simple triangulations, let us 
conclude this section by expressing $R_\infty$, $S_\infty$ and $g$ in terms of the parameter $\lambda$.
Introducting the intermediate variable $s\equiv S_\infty/\sqrt{R_\infty}$, we may write \eqref{eq:recRS} as
\begin{equation*}
\left\{
\begin{split}
& R_\infty =1+2 \sqrt{g^2\, R_\infty^3}\,  s\\
& R_\infty\, s = \sqrt{g^2\, R_\infty^3} (s^2 +2)\\
\end{split}
\right.
\end{equation*}
which, after eliminating $s$ from the system, implies
\begin{equation*}
 R_\infty^2=1+8\, (g^2\, R_\infty^3)\ .
\end{equation*}
From \eqref{eq:charac}, this leads to
\begin{equation*}
R_\infty=\left(1+ \frac{8}{\lambda+\frac{1}{\lambda}+2}\right)^{1/2}=\frac{\sqrt{1+10\lambda+\lambda^2}}{1+\lambda}
\end{equation*}
with $\lambda$ parametrizing $g$ via
\begin{equation*}
g= \frac{1}{R_\infty^{3/2}\, \left(\lambda+\frac{1}{\lambda}+2\right)^{1/2}}=\frac{\sqrt{\lambda\, (1+\lambda)}}{(1+10\lambda+\lambda^2)^{3/4}}\ .
\end{equation*}
Note that, since $\lambda\neq 1$, $g$ ranges from $0$ to $1/(2\cdot 3^{3/4})$. in agreement with the fact that the number of triangulations 
with $n$ faces growths like $(2\cdot 3^{3/4})^n$ (see for instance \cite{Gao91} for explicit formulas).
Finally we have
\begin{equation*}
S_\infty=\frac{R_\infty-1}{2g R_\infty}=(1+10\lambda+\lambda^2)^{1/4}\, \frac{\sqrt{1+10\lambda+\lambda^2}-(1+\lambda)}{2\sqrt{\lambda\, (1+\lambda)}}\ .
\end{equation*}
Plugging the above formulas for $R_k$ and $S_k$ in eq.~\eqref{eq:twopoint} and the expressions for $g$, $R_\infty$
and $S_\infty$, we arrive at the remarkably simple expression of the distance-dependent two-point function
$G_k$ for $k\geq 1$:
\begin{equation*}
G_k=\frac{(1-\lambda^3)(1+10\, \lambda+\lambda^2)}{1+\lambda}\, 
\frac{\lambda^{k-1}(1+\lambda^{k+1})}{(1-\lambda^k)(1-\lambda^{k+1})(1-\lambda^{k+2})}-\delta_{k,1}\ .
\end{equation*}

\section{Conclusion}
\label{sec:conclusion}
In this paper, we presented a new technique to compute the distance-dependent two-point function of planar
triangulations by first deriving and then solving a new recursion for the intimately related slice generating functions.
Although our method makes a crucial use of properties which are specific to triangulations, it is likely that it could be
generalized to other families of maps. In particular, the case of \emph{planar quadrangulations} seems promising for
a similar treatment.

Our approach is based on the existence, in slices of left-boundary length $\ell$, of
a dividing line connecting the right and left boundaries of the slice via $\ell-1 \to \ell-1$ edges.
Upon gluing, say the two boundaries of an $R$-slice with $\ell=k$, we produce via the equivalence displayed 
in figure~\ref{fig:twopoint1} a pointed rooted triangulation whose root-edge is ``of type" $k\to k-1$ 
with respect to the marked vertex. After gluing, the dividing line creates a simple closed path made of edges 
connecting vertices at distance $k-1$ from the marked vertex, and which separates the 
marked vertex from the root-vertex. By construction, all the vertices strictly outside the domain containing the 
marked vertex are at distance at least $k$ from this vertex. We may thus interpret the dividing line as the
boundary of the \emph{hull of radius $k-1$} centered at the marked vertex which, so to say, is formed of the
ball of radius $k-1$ together with \emph{all the complementary connected domains} (thus containing vertices at 
distance at least $k$ 
from the marked vertex) \emph{except that containing the root-vertex}. In other words, we may decide to use 
the dividing line as a way to precisely \emph{define}
what we shall call the hull of radius $k-1$ centered at the marked vertex, namely the domain delimited by this line and containing
the marked vertex. Note that this definition is \emph{not equivalent} to that given in \cite{Krikun03} where the hull 
of radius $k-1$ is constructed explicitly from the ball of radius $k-1$, itself defined at the set of all triangles incident 
to at least one vertex at distance $\leq k-2$ from the marked vertex (see \cite{Krikun03}). 

Finally, our recursive construction allows us to define dividing lines within each of the slices (enumerated by $T_{k-1}$)
which appear in the slice decomposition of the hull (see figure~\ref{fig:twopoint1}). After gluing, the concatenation 
of these dividing lines creates a simple closed path made of edges 
connecting vertices at distance $k-2$ from the marked vertex, and which separates the 
marked edge from root-vertex. This line may be viewed as the boundary of the hull of radius $k-2$ centered at 
the marked vertex.
We may in this way define hulls of all radii between $1$ and $k-1$ and our recursion relations should
in principle allow us to describe the statistics of the lengths of these hull boundaries, an analysis yet to be done.



\bibliographystyle{plain}
\bibliography{triang2p}

\begin{thebibliography}{10}

\bibitem{AmBudd}
J.~{Ambj\o rn} and T.G. Budd.
\newblock Trees and spatial topology change in causal dynamical triangulations.
\newblock {\em J. Phys. A: Math. Theor.}, 46(31):315201, 2013.

\bibitem{GEOD}
J.~Bouttier, P.~Di~Francesco, and E.~Guitter.
\newblock Geodesic distance in planar graphs.
\newblock {\em Nucl. Phys. B}, 663(3):535--567, 2003.

\bibitem{BFG}
J.~Bouttier, \'E. Fusy, and E.~Guitter.
\newblock On the two-point function of general planar maps and hypermaps.
\newblock {\em Ann. Inst. Henri Poincar\'e Comb. Phys. Interact.},
  1(3):265--306, 2014.
\newblock arXiv:1312.0502 [math.CO].

\bibitem{BG12}
J.~Bouttier and E.~Guitter.
\newblock Planar maps and continued fractions.
\newblock {\em Comm. Math. Phys.}, 309(3):623--662, 2012.

\bibitem{CoriVa}
R.~Cori and B.~Vauquelin.
\newblock Planar maps are well labeled trees.
\newblock {\em Canad. J. Math.}, 33(5):1023--1042, 1981.

\bibitem{PDFRaman}
P.~Di Francesco.
\newblock Geodesic distance in planar graphs: an integrable approach.
\newblock {\em Ramanujan J.}, 10(2):153--186, 2005.

\bibitem{Gao91}
Z.-C. Gao.
\newblock The number of rooted triangular maps on a surface.
\newblock {\em Journal of Combinatorial Theory, Series B}, 52(2):236--249,
  1991.

\bibitem{Krikun03}
M.A. Krikun.
\newblock Uniform infinite planar triangulation and related time-reversed
  critical branching process.
\newblock {\em Journal of Mathematical Sciences}, 131(2):5520--5537, 2005.

\bibitem{SchPhD}
G.~Schaeffer.
\newblock {\em Conjugaison d'arbres et cartes combinatoires al\'eatoires}.
\newblock PhD thesis, Universit\'e Bordeaux I, 1998.

\bibitem{TutteCPT}
W.~T. Tutte.
\newblock A census of planar triangulations.
\newblock {\em Canad. J. Math.}, 14:21--38, 1962.

\end{thebibliography}

\end{document}